\documentclass{article}

\usepackage[utf8]{inputenc}
\usepackage{amsmath}
\usepackage{mathtools}
\usepackage{amsfonts}
\usepackage{amssymb}
\usepackage{amsthm}

\usepackage{hyperref}

\usepackage[compress,sort]{cite}

\usepackage{cases}

\newcommand{\dlsigma}{\delta\sigma}
\newcommand{\dlv}{\delta v}
\newcommand{\dlw}{\delta w}

\newcommand{\dt}{\partial_t}
\newcommand{\nablah}{\nabla_h}
\newcommand{\deltah}{\Delta_h}

\newcommand{\ddh}{\partial_h}
\newcommand{\dz}{\partial_z}

\newcommand{\idx}{\,d\vec{x}}
\newcommand{\idxh}{\,d\vec{x}_h}

\newcommand{\dvh}{\mathrm{div}_h}

\newcommand{\subeqref}[2]{$ \eqref{#1}_{#2} $}
\newcommand{\norm}[2]{\Arrowvert #1 \Arrowvert_{#2}}
\newcommand{\hnorm}[2]{| #1 |_{#2}}
\newcommand{\Lnorm}[1]{L^{#1}}
\newcommand{\Hnorm}[1]{H^{#1}}

\theoremstyle{theorem}

\newtheorem{thm}{Theorem}

\newtheorem{prop}{Proposition}

\theoremstyle{remark}
\newtheorem{rmk}{Remark}

\allowdisplaybreaks

\title{Rigorous justification of the hydrostatic approximation limit of viscous compressible flows}

\author{Xin Liu\footnote{Department of Mathematics, Texas A{\&}M University, College Station, TX 77843, USA. Email: xliu23@tamu.edu} \,\, and \,  Edriss S. Titi\footnote{Department of Mathematics, Texas A{\&}M University, College Station,  TX 77843, USA.  Department of Applied Mathematics and Theoretical Physics, University of Cambridge, Cambridge CB3 0WA, UK.
		Department of Computer Science and Applied Mathematics, Weizmann Institute of Science, Rehovot 76100, Israel. Email: {titi@math.tamu.edu}\, and \, {Edriss.Titi@maths.cam.ac.uk}}
}

\date{\today}

\begin{document}
	\maketitle

\begin{abstract}
	This paper considers the asymptotic limit of small aspect ratio between vertical and horizontal spatial scales for viscous isothermal compressible flows. 
	In particular, it is observed that fast vertical acoustic waves arise and induce an averaging mechanism of the density in the vertical variable, which at the limit leads to the hydrostatic approximation of compressible flows, i.e., the compressible primitive equations of atmospheric dynamics. We justify the hydrostatic approximation for general as well as ``well-prepared'' initial data. The initial data is called well-prepared when it is close to the hydrostatic balance in a strong topology. Moreover, the convergence rate is calculated in the well-prepared initial data case in terms of the aspect ratio, as the latter goes to zero. 
	
\end{abstract}

\tableofcontents
	
\section{Introduction}

	\subsection{Hydrostatic approximation of compressible flows}
		
	The compressible primitive equations of atmospheric dynamics are the hydrostatic approximation of the compressible hydrodynamic equations, which are obtained by replacing the evolutionary vertical momentum equation with the hydrostatic balance equation (see \subeqref{CPE}{3}, below). A formal derivation can be found, e.g., in \cite{Ersoy2011a}. The small aspect ratio in the atmosphere between the vertical scale and the horizontal planetary scale plays an essential role in this derivation and is the main factor behind the hydrostatic approximation. The hydrostatic approximation is commonly used in the atmospheric science models, and has successfully simplified the complex hydrodynamic equations and their computational aspects. In this work, our goal is to rigorously justify the hydrostatic approximation for viscous compressible flows and to provide a mathematical foundation for its application. For more backgrounds about the hydrostatic approximation, see, e.g., \cite{LT2018a,LT2018b,Lions1992,JLLions1994}
	
	Recall that the dynamics of compressible isothermal flow in a thin periodic channel domain is governed, after omitting the viscosities, by the isothermal compressible Euler equations in the domain $ \Omega_\varepsilon = 2\mathbb T^2 \times (0, \varepsilon) $, for small $ \varepsilon \in (0,1) $:
	\begin{equation*}\tag{EQ}\label{EQ}
		\begin{cases}
			\dt \rho + \dvh (\rho v) + \dz (\rho w) = 0, & \text{in} ~  \Omega_\varepsilon, \\
			\rho \dt v + \rho v\cdot \nablah v + \rho w \dz v + \nablah \rho = 0, & \text{in} ~ \Omega_\varepsilon, \\
			\rho \dt w + \rho v\cdot \nablah w + \rho w \dz w + \dz \rho = 0, & \text{in} ~ \Omega_\varepsilon,
		\end{cases}
	\end{equation*}
	subject to the impenetrable boundary condition at the wall boundaries
	\begin{equation*}
		w\big|_{z = 0, \varepsilon} = 0,
	\end{equation*}
	where $ \rho, v, w $ represent the unknown density, the horizontal velocity, and the vertical velocity, respectively. Here $ \dvh, \nablah $ represent the horizontal divergence and the horizontal gradient, respectively (see section \ref{subsec:preliminary}, below). $ 2 \mathbb T^2 \subset \mathbb R^2 $ is the periodic horizontal domain (flat torus) with period $ 2 $ in each direction. 
	
	Alternatively, set $ \sigma  : = \log \rho $, which, for simplicity, we will still refer to as the density. Then, in the region where $ \rho > 0 $, \eqref{EQ} is equivalent to
	\begin{equation*}\tag{\ref{EQ}'}\label{EQ'}
		\begin{cases}
			\dt \sigma + v \cdot \nablah \sigma + w \dz \sigma + \dvh v + \dz w = 0, & \text{in} ~ \Omega_\varepsilon, \\
			\dt v + v \cdot \nablah v + w \dz v + \nablah \sigma = 0, & \text{in} ~ \Omega_\varepsilon, \\
			\dt w + v \cdot \nablah w + w \dz w + \dz \sigma = 0, & \text{in} ~ \Omega_\varepsilon, 
		\end{cases}
	\end{equation*}
	with
		$ w \big|_{z=0,\varepsilon} = 0. $
	As $ \varepsilon \rightarrow 0^+ $, it is expected that $$ ( \dfrac{1}{\varepsilon}\int_0^\varepsilon \sigma(\cdot,z')\,dz', \dfrac{1}{\varepsilon}\int_0^\varepsilon v(\cdot, z')\,dz') $$ will converge to solutions to the two-dimensional compressible Euler equations. This has been verified in the viscous case (i.e., the compressible Navier-Stokes equations) in \cite{Bella2014} (see also \cite{Maltese2014}). In view of multi-scale analysis, this is to say that, if one considers the multi-scale expansion of the solutions to \eqref{EQ'} by writing
	\begin{align*}
		\sigma(x,y,z,t) = & \sigma^0(x,y,z/\varepsilon,t) + \varepsilon \sigma^1(x,y,z/\varepsilon,t) + \cdots, \\
		v(x,y,z,t) = & v^0(x,y,z/\varepsilon,t) + \varepsilon v^1(x,y,z/\varepsilon,t) + \cdots, \\
		w(x,y,z,t) = & w^0(x,y,z/\varepsilon,t) + \varepsilon w^1(x,y,z/\varepsilon,t) + \cdots, 
	\end{align*}
	it should be expected that $ ( \int_0^1 \sigma^0(\cdot, z')\,dz', \int_0^1 v^0(\cdot,z')\,dz') $ solves the two-dimensional compressible Euler equations and $ w^0 \equiv 0 $. Therefore, the non-trivial leading order of the vertical velocity $ w $ is $ w^1 $, and one can easily check that the equations satisfied by $ (\sigma^0, v^0, w^1) $ are exactly given by the inviscid compressible primitive equations, i.e., equations \eqref{CPE}, below, without the viscosities. To capture the above scale analysis rigorously, we consider the ansatz
	\begin{equation}\label{ansatz}
		\begin{aligned}
			\sigma(x,y,z) : = & \sigma_\varepsilon(x,y,z/\varepsilon), \\
			v(x,y,z) := & v_{\varepsilon} (x,y,z/\varepsilon), \\
			w(x,y,z) := & \varepsilon w_{\varepsilon} (x,y,z/\varepsilon),
		\end{aligned}
	\end{equation}
	and denote by $ z' := z/\varepsilon \in (0,1) $ and $ \Omega := 2\mathbb T^2 \times (0,1) $. Then we write down the equations satisfied by $ (\sigma_\varepsilon, v_\varepsilon, w_\varepsilon) $ in $ \Omega $, taking into account the eddy viscosities:
		\begin{equation}\label{CF-0}
			\begin{cases}
			\dt \sigma_{\varepsilon} + v_{\varepsilon}\cdot\nablah \sigma_{\varepsilon} + w_{\varepsilon} \dz \sigma_{\varepsilon} + \dvh v_{\varepsilon} + \dz w_{\varepsilon} = 0, & \text{in} ~ \Omega, \\
			\dt v_{\varepsilon} + v_{\varepsilon} \cdot \nablah v_{\varepsilon} + w_{\varepsilon} \dz v_{\varepsilon} + \nablah \sigma_{\varepsilon} = \deltah v_{\varepsilon} + \partial_{zz} v_{\varepsilon}, & \text{in} ~ \Omega,\\
			\varepsilon^2 \bigl(\dt w_{\varepsilon} + v_{\varepsilon} \cdot \nablah w_{\varepsilon} + w_{\varepsilon} \dz w_{\varepsilon} \bigr) + \dz \sigma_{\varepsilon} = \varepsilon^2 \bigl( \deltah w_{\varepsilon} + \partial_{zz} w_{\varepsilon} \bigr), & \text{in} ~ \Omega,
			\end{cases}
		\end{equation}
		subject to the impenetrable and stress-free boundary conditions at the wall boundaries
		\begin{equation}\label{bc-CF-0}
			w_\varepsilon\big|_{z=0,1} = 0, ~ \dz v_\varepsilon\big|_{z=0,1} = 0,
		\end{equation}
		where we have dropped the prime sign for the vertical variable. 
		
		We remark that \eqref{CF-0} can also be obtained by considering ansatz \eqref{ansatz} in the viscous version of \eqref{EQ'} with viscosities $ \deltah v_\varepsilon + \varepsilon^{-2} \partial_{zz} v_\varepsilon $ and $ \deltah w_\varepsilon + \varepsilon^{-2} \partial_{zz} w_\varepsilon $ on the right-hand side of the horizontal momentum equation \subeqref{EQ'}{2} and the vertical momentum equation \subeqref{EQ'}{3}, respectively. 
		
		Let us observe that the solutions to \eqref{CF-0} is invariant with respect to the following symmetry:
		\begin{equation*}\tag{SYM}\label{SYM}
			\text{$ \sigma_\varepsilon, v_\varepsilon, w_\varepsilon $ are even, even, and odd in the $ z $-variable, respectively.}
		\end{equation*}
		For this reason, in this work, we equivalently consider the following viscous compressible hydrodynamic system with turbulence eddy viscosities: for small $ \varepsilon \in (0,1) $,
		\begin{equation}\label{CF}
			\begin{cases}
			\dt \sigma_{\varepsilon} + v_{\varepsilon}\cdot\nablah \sigma_{\varepsilon} + w_{\varepsilon} \dz \sigma_{\varepsilon} + \dvh v_{\varepsilon} + \dz w_{\varepsilon} = 0, & \text{in} ~ 2 \mathbb T^3, \\
			\dt v_{\varepsilon} + v_{\varepsilon} \cdot \nablah v_{\varepsilon} + w_{\varepsilon} \dz v_{\varepsilon} + \nablah \sigma_{\varepsilon} = \deltah v_{\varepsilon} + \partial_{zz} v_{\varepsilon}, & \text{in} ~ 2 \mathbb T^3,\\
			\varepsilon^2 \bigl(\dt w_{\varepsilon} + v_{\varepsilon} \cdot \nablah w_{\varepsilon} + w_{\varepsilon} \dz w_{\varepsilon} \bigr) + \dz \sigma_{\varepsilon} = \varepsilon^2 \bigl( \deltah w_{\varepsilon} + \partial_{zz} w_{\varepsilon} \bigr), & \text{in} ~ 2 \mathbb T^3,
			\end{cases}
		\end{equation}
		satisfying the symmetry in \eqref{SYM} and the periodic boundary condition in all directions. Notice that the boundary conditions in \eqref{bc-CF-0} are automatically satisfied by regular solutions to \eqref{CF} owing to symmetry \eqref{SYM}, and one can obtain a solution to \eqref{CF-0} with \eqref{bc-CF-0} by restricting any regular solution of \eqref{CF} to the domain $ \Omega $. 
		
		The formal limit system of \eqref{CF}, as $ \varepsilon \rightarrow 0^+ $, is
		\begin{equation}\label{CPE}
			\begin{cases}
			\dt \sigma_p + v_p \cdot\nablah \sigma_p + w_p \dz \sigma_p + \dvh v_p + \dz w_p = 0, & \text{in} ~ 2\mathbb T^3, \\
			\dt v_p + v_p \cdot \nablah v_p + w_p \dz v_p + \nablah \sigma_p = \deltah v_p + \partial_{zz} v_p, & \text{in} ~ 2\mathbb T^3,\\
			\dz \sigma_p = 0, & \text{in} ~ 2\mathbb T^3,
			\end{cases}
		\end{equation}
		satisfying the same symmetry as in \eqref{SYM}, with $ \sigma_\varepsilon, v_\varepsilon, w_\varepsilon $ replaced by $ \sigma_p, v_p, w_p $, respectively. 

\bigskip
		
		Before we move on to the discussion of our strategy to rigorously establish the aforementioned limiting problem, we would like to mention some relevant previous results. 
		
		As for the compressible hydrodynamic equations for viscous flows, i.e., the compressible Navier-Stokes equations, the derivation of the system can be found in \cite{Lions1996,Feireisl2004}. In \cite{Lions1998,Feireisl2004}, global weak solutions to the compressible Navier-Stokes equations are constructed. Recently, the authors in \cite{Li2015a} and \cite{Vasseur2016}  independently construct global weak solutions to the compressible Navier-Strokes equations with degenerate viscosities. See also \cite{Bresch2006} and \cite{Bresch2019} for relevant developments. As for the strong solutions of the compressible Navier-Stokes equations, \cite{Choe2003,Cho2004,Cho2006a,Cho2006c} establish the local well-posedness of strong and classical solutions with vacuum. Without vacuum, the local well-posedness theory can be dated back to \cite{Serrin1959,Itaya1971,Tani1977}. The first global well-posedness result is established in \cite{Matsumura1980,Matsumura1983}, where the asymptotic stability of constant states is studied with respect to small perturbations. A global existence theorem with vacuum and small energy is given in \cite{HuangLiXin2012} (see also \cite{Huang2018a}). We refer readers, for other developments of the compressible Navier-Stokes equations to, e.g.,  \cite{Vaygant1995,zpxin1998,Huang2016,Hoff1992,hoff1991,Bresch2018}.
		
		As for the compressible primitive equations (PE) \eqref{CPE}, one can find the meteorological discussion and applications of the system in, e.g., \cite{Richardson1965} and \cite{Washington2005}. In two dimensions, the global weak solutions are constructed in \cite{Gatapov2005,Ersoy2012}. The stability of weak solutions is established in \cite{Ersoy2011a}. Uniqueness of the weak solutions in two dimensions is studied in \cite{Jiu2018}. As for the three-dimensional dynamics, the existence of global weak solutions for the compressible primitive equations with degenerate viscosities is established in \cite{LT2018b} and \cite{Wang2017}, independently. The local well-posedness of strong solutions for the compressible primitive equations with constant viscosities is established in \cite{LT2018a}. We introduce the PE diagram and study the small Mach number limit of the compressible primitive equations in \cite{LT2018LowMach1} and \cite{LT2022LowMath2}.  For readers' convenience, results concerning the incompressible primitive equations can be found in \cite{HuTemamZiane2003,Petcu2005,GuillenGonzalez2001,Kukavica2007a,Cao2007,Renardy2009,Cao2012,Cao2014b,Cao2014a,Cao2016,Cao2017,Cao2016a,Li2017a,GerardVaret2018,Li2017,Azerad2001}. 
		
		We would like to remark that, the well-posedness of local strong or classical solutions to \eqref{CF} for any fixed $ \varepsilon \in (0,1) $, and \eqref{CPE}, follows easily, using the similar arguments as in \cite{Itaya1971} and \cite{LT2018a}, respectively. We only mention that $ w_p $ in \eqref{CPE} can be calculated by 
		\begin{equation}\label{id:vertical-p}
			w_p(\cdot_h,z) = - e^{-\sigma_p(\cdot_h)} \int_0^z e^{\sigma_p(\cdot_h)} \biggl( \widetilde v_p(\cdot_h,z') \cdot \nablah \sigma_p(\cdot_h) + \dvh \widetilde v_p (\cdot_h,z')\biggr)\,dz' ,
		\end{equation}
		using the continuity equation \subeqref{CPE}{1}, while the vertical average part of the continuity equation yields an evolutionary equation for $ \sigma_p $, i.e., 
		\begin{equation*}
			\dt \sigma_p + \overline v_p\cdot\nablah \sigma_p + \dvh \overline v_p = 0. 
		\end{equation*}
		Here we have used the fact that $ \sigma_p $ is independent of the $ z $-variable due to the hydrostatic equation \subeqref{CPE}{3}, $\cdot_h$ denotes the horizontal variables, i.e., $x,y$, and $\overline \cdot $ and $ \widetilde \cdot $ are the vertical average and fluctuation, respectively, defined explicitly in \eqref{def:avg-flc} in section \ref{subsec:preliminary}, below. 

\bigskip
		
		To establish the asymptotic limit as $ \varepsilon \rightarrow 0^+ $, we will need to establish some uniform-in-$\varepsilon $ estimates for $ (\sigma_\varepsilon, v_\varepsilon, w_\varepsilon) $. While in the case of the incompressible flows, e.g., \cite{Li2017,Azerad2001,Li-Titi-Yuan2022}, the vertical velocity can be represented by the horizontal velocity, using the incompressibility condition, this benefit no longer exists in the compressible case. Instead, for the compressible primitive equations, as in \eqref{id:vertical-p}, the vertical velocity is represented by the density $ \sigma_p $ combined with the horizontal velocity $ v_p $, which does not involve any time derivative. This allows us to consider \subeqref{CPE}{1} and \subeqref{CPE}{2} as closed evolutionary equations of $ (\sigma_p, v_p) $ only, and thus allows us to construct the local strong solutions in \cite{LT2018a}. However, all these structures do not work for equations \eqref{CF}. Indeed, the vertical momentum equation \subeqref{CF}{3} can only provide estimates for $ \varepsilon w_\varepsilon $, while using the continuity equation \subeqref{CF}{1}, the representation of $ w_\varepsilon $ in terms of $ \sigma_\varepsilon $ and $ v_\varepsilon $ involves the time derivative $ \dt \sigma_\varepsilon $ (see \eqref{id:vertical-velocity-0}, below). This may involve high oscillations, as $ \varepsilon \rightarrow 0^+ $, and is the main obstacle to overcome in this work. It is worth noticing that, in \cite{gao2022hydrostatic,tang2023derivation}, by assuming the uniform-in-$\varepsilon$ existence of, and bounds on, global weak solutions to fluid system \eqref{CF-0} with and without viscosity, and with the additional assumption that the density is independent of the $ z $-variable, the authors justify the hydrostatic approximation in compressible fluid. However, these assumptions in \cite{gao2022hydrostatic,tang2023derivation} can not be justified physically for atmospheric dynamics. They are mathematically assumed in order to essentially avoid dealing with the obstacle mentioned above. In this paper, our goal is to overcome these difficulties by showing the uniform-in-$\varepsilon$ existence, bounds, and convergence. 
		
		\bigskip
		To motivate our strategy and to shed light on the treatment of this obstacle, let us consider the following linear model system:
		\begin{equation}\label{eq:linear-model}
			\begin{cases}
				\dt \eta + \dvh \psi^h + \dz \psi^z =0, \\
				\dt \psi^h + \nablah \eta = \Delta \psi^h , \\
				\varepsilon^2 \dt \psi^z + \dz \eta = \varepsilon^2 \Delta \psi^z,
			\end{cases}
		\end{equation}
		with unknown scalar functions $ \eta, \psi^h$, and $ \psi^z $ in $ 2 \mathbb T^3 $, subject to the symmetry that $ \eta, \psi^h, \psi^z $ are even, even, and odd in the $ z $-variable, respectively, similar to \eqref{SYM}. In particular, $ \psi^z|_{z=0} = 0 $. 
		Then the standard $ H^s $ estimate of \eqref{eq:linear-model}, for every $ s \in \mathbb Z^+ $, yields
		\begin{equation*}
			\sup_{0\leq t\leq T}  \norm{\eta(t),\psi^h(t), \varepsilon\psi^z(t)}{\Hnorm{s}}^2 + \int_0^T \norm{\psi^h(t),\varepsilon\psi^z(t)}{\Hnorm{s+1}}^2 \,dt < \infty,
		\end{equation*}
		for any $ T \in (0,\infty) $, which allows us to pass the limit $ \varepsilon \rightarrow 0^+ $ in \subeqref{eq:linear-model}{2} and \subeqref{eq:linear-model}{3}, but not \subeqref{eq:linear-model}{1}, due to the lack of compactness of the $ \psi^z $ sequence. 

\bigskip		

		To overcome this obstacle, we focus instead on the hyperbolic structure of the system
		\begin{equation}\label{eq:linear-model-2}
			\begin{cases}
				\dt \eta + \dz \psi^z  +\cdots = 0, \\
				\dt \psi^z + \dfrac{1}{\varepsilon^2} \dz \eta + \cdots = 0. 
			\end{cases}
		\end{equation}
		Then from \eqref{eq:linear-model-2}, one can obtain a wave equation for $ \eta $, i.e., 
		\begin{equation*}
			\partial_{tt} \eta - \dfrac{1}{\varepsilon^2} \partial_{zz} \eta = \cdots,
		\end{equation*}
		which can provide uniform-in-$\varepsilon $ estimates for $ \dt \eta $. In the end, using \subeqref{eq:linear-model-2}{1} again, one can write $ \psi^z(z)  = - \int_0^z ( \dt \eta(z') + \cdots ) \,dz' $, from which one can obtain the required uniform-in-$\varepsilon$ estimates of $ \psi^z $, and thus the missing (weak) compactness of $ \psi^z $ is obtained. In other words, one has to take advantage of the oscillatory nature of the underlying system in order to obtain the required uniform-in-$ \varepsilon $ estimates.
		
\bigskip

		Back to \eqref{eq:linear-model}, the contribution of the viscosities (i.e., $\Delta \psi^h$ and $ \Delta \psi^z $) should also be taken into consideration. One can calculate, similarly, 
		\begin{equation*}
			\partial_{tt} \eta = - \dvh \dt \psi^h - \dz \dt \psi^z = \deltah \eta + \dfrac{1}{\varepsilon^2} \partial_{zz}\eta - \Delta (\dvh \psi^h + \dz \psi^z),
		\end{equation*}
		where, using \subeqref{eq:linear-model}{1}, we have $ \Delta (\dvh \psi^h + \dz \psi^z) = - \dt \Delta \eta $. Therefore, we arrive at
		\begin{equation*}
			\dt (\dt \eta - \Delta \eta) - \deltah \eta - \dfrac{1}{\varepsilon^2}\partial_{zz}\eta = 0,
		\end{equation*}
		which is a strongly damped wave equation. Using such a structure in the nonlinear problem \eqref{CF} (see \eqref{eq:mixed}, below), one can obtain the required uniform-in-$\varepsilon$ estimates for $ \dt \sigma_\varepsilon $ and hence those for $ w_\varepsilon $. 
		
		In order to deal with the nonlinearities, we complement \eqref{CF} with $ H^3 $ initial data. However, we remark that this may not be the optimal regularity for the initial data. 

\bigskip
		
		After obtaining the aforementioned uniform-in-$\varepsilon$ estimates, we will be able to establish the limit of \eqref{CF} as $ \varepsilon \rightarrow 0^+ $. This is done in section \ref{sec:convergence-cpe}. However, this is not enough to establish the convergence rates as $ \varepsilon \rightarrow 0^+ $. To explain why, consider the difference 
$ (\dlsigma, \dlv, \dlw) := (\sigma_{\varepsilon} - \sigma_p, v_{\varepsilon}-v_p, w_{\varepsilon}-w_p) $. Then $ (\dlsigma, \dlv, \dlw) $ satisfies
		\begin{equation}\label{eq:difference}
			\begin{cases}
			\dt \dlsigma + v_{\varepsilon} \cdot \nablah \dlsigma + w_{\varepsilon} \dz \dlsigma + \dlv \cdot \nablah \sigma_p + \dvh \dlv + \dz \dlw = 0,\\
			\dt \dlv + v_{\varepsilon} \cdot \nablah \dlv + w_{\varepsilon} \dz \dlv + \dlv \cdot \nablah v_p + \dlw \dz v_p \\
			~~~~ ~~~~ + \nablah \dlsigma = \deltah \dlv + \partial_{zz} \dlv, \\
			\dt \dlw + v_{\varepsilon} \cdot\nablah \dlw + w_{\varepsilon} \dz \dlw + \dlv \cdot \nablah w_p + \dlw \dz w_p \\
			~~~~ ~~~~ + \dfrac{1}{\varepsilon^2} \dz \dlsigma = \deltah \dlw + \partial_{zz} \dlw 
			+ \bigl( \deltah w_p + \partial_{zz} w_p - \dt w_p\\
			~~~~ ~~~~ - v_p\cdot \nablah w_p - w_p \dz w_p  \bigr).
			\end{cases}
		\end{equation}
		While the uniform-in-$\varepsilon $ estimates work well for $ (\dlsigma, \dlv, \dlw ) $, because of the terms $$ \deltah w_p + \partial_{zz} w_p - \dt w_p - v_p\cdot\nablah w_p - w_p\dz w_p $$
		on the right-hand side of \subeqref{eq:difference}{3}, the uniform-in-$\varepsilon$ estimates are not comparable to $ \varepsilon $, with the exception, however, only of the estimate 
		\begin{equation}\label{fluctuation-est}
			\dz \dlsigma \sim \mathcal O(\varepsilon) ~~ \text{in some norm.}
		\end{equation}
		See \eqref{uniform-est-2}, below. For comparison with the case of incompressible flows, we refer the reader to \cite{Li2017,Li-Titi-Yuan2022}, where such an issue does not exist. 
		
		Writing $ \dlsigma = \overline{\dlsigma} + \widetilde{\dlsigma} $, see \eqref{def:avg-flc}, below, we notice that \eqref{fluctuation-est} implies that the fluctuation of $ \dlsigma $ is of order $ \varepsilon$, i.e., $ \widetilde \dlsigma \sim \mathcal O(\varepsilon) $ in some norm. Inspired by the study of \eqref{CPE} in \cite{LT2018a}, one can separate \subeqref{eq:difference}{1} into the (vertical) average part and the fluctuation part. Then the average part is nothing but an evolutionary equation of $ \overline \dlsigma $ and does not involve $ \dlw $. Therefore, using the average part and \subeqref{eq:difference}{2}, one can show that $ \overline \dlsigma \sim \mathcal O(\varepsilon) $ and $ \dlv \sim \mathcal O(\varepsilon) $ in some norm. On the other hand, the fluctuation part yields that the estimates for $ \dlw $, in fact,  only involve $ \dt \widetilde \dlsigma $. Therefore, one will only need to obtain an estimate of the order $ \dz\dt\dlsigma \sim \mathcal O(\varepsilon) $ in some norm in order to close the estimates of converging rates. This is established by the use of some additional uniform-in-$\varepsilon $ estimates in section \ref{subsec:uniform-est-2}. That is, by assuming a well-prepared initial data, i.e., initial data that is close to the hydrostatic approximation (see assumption \eqref{WP}, below).

\bigskip
		
		The rest of this paper is organized as follows. In section \ref{subsec:preliminary}, we introduce the notations that have been and will be used in this paper, as well as the energy and dissipation functionals. In section \ref{subsec:theorem}, we summarize the main theorems. In section \ref{sec:uniform-est}, we establish the uniform-in-$\varepsilon $ estimates, which imply the necessary weak and strong compactness in order to pass the limit $ \varepsilon \rightarrow 0^+ $ in section \ref{sec:convergence-cpe}. In section \ref{sec:converging-rate}, we focus on the study of convergence rates, which is established for a restricted class of well-prepared initial data, see \eqref{WP}, below. Readers who are more interested in the converging rates can skip directly to section \ref{subsec:converging-rate}.

	\subsection{Preliminaries}\label{subsec:preliminary}
	In this paper, we consistently use $ t  \in [0,\infty) $ to represent the temporal variable, $ x, y \in 2\mathbb T $ to represent the horizontal spatial variables, and $ z \in 2 \mathbb T $ to represent the vertical spatial variable. 
	We have and will use $ \nabla_h, \dvh $, and $ \Delta_h $ to represent the horizontal gradient, the horizontal divergence, and the horizontal Laplace operator, respectively; that is,
	\begin{gather*}
		\nabla_h := \biggl(\begin{array}{c}
			\partial_x \\ \partial_y
		\end{array}\biggr) , ~ \dvh := \nabla_h \cdot, ~
		\Delta_h := \dvh \nabla_h.
	\end{gather*}
	Using such notations, we have
	$$ \Delta = \deltah + \partial_{zz}. $$
	Also, for any function $ f $, we denote the vertical average and fluctuation of $ f $ as 
	\begin{equation}\label{def:avg-flc}
		\overline{f}(x,y):= \int_0^1 f(x,y,z)\,dz, ~ \widetilde{f} := f- \overline{f},
	\end{equation}
	also known as the barotropic and baroclinic modes, respectively. 
	$ \hnorm{f(z,t)}{X}, \norm{f(t)}{X} $ are used to denote the $X$-norm in the horizontal domain $  2 \mathbb T^2  $ for any fixed $ z, t $, and in the three-dimensional domain $ 2\mathbb T^3 $ for any fixed $ t $, respectively.

The following functionals will be used in this paper:
	\begin{align}
	\label{energy}
	& \begin{aligned} 
	& E = E(t) := \norm{v_{\varepsilon}, \varepsilon w_{\varepsilon}}{\Hnorm{3}} + \norm{\dt \sigma_{\varepsilon},\nablah \sigma_{\varepsilon}, \dfrac{\dz \sigma_{\varepsilon}}{\varepsilon}}{\Hnorm{2}} + \norm{\sigma_{\varepsilon}}{\Hnorm{4}} \\
	& ~~~~ ~~~~ + \norm{w_{\varepsilon}, \dz w_{\varepsilon}}{\Hnorm{2}} , \end{aligned} \\
	\label{dissipation}
	& \begin{aligned}
	&  D = D(t) := \norm{\dt v_{\varepsilon},  \dt (\varepsilon w_{\varepsilon})}{\Hnorm{2}} + \norm{v_{\varepsilon}, \varepsilon w_{\varepsilon}}{\Hnorm{4}} + \norm{\dt \sigma_{\varepsilon}, \nablah \sigma_{\varepsilon}, \dfrac{ \dz \sigma_{\varepsilon}}{\varepsilon} }{\Hnorm{3}} \\
	& ~~~~ ~~~~  + \norm{w_{\varepsilon}, \dz w_{\varepsilon}}{\Hnorm{3}}.  
	\end{aligned}
	\end{align}
	Correspondingly, we denote the functionals of temporal derivatives as,
	\begin{align}
	\label{energy-1}
	& \begin{aligned} 
	& E_1 = E_1(t) := \norm{\dt v_{\varepsilon}, \dt(\varepsilon w_{\varepsilon})}{\Hnorm{1}} + \norm{\dt^2 \sigma_{\varepsilon}, \dt \nablah \sigma_{\varepsilon}, \dfrac{\dz \dt \sigma_{\varepsilon}}{\varepsilon}}{\Lnorm{2}} \\
	& ~~~~ ~~~~ + \norm{\dt \sigma_{\varepsilon}}{\Hnorm{2}} + \norm{\dt w_{\varepsilon}, \dz \dt w_{\varepsilon} }{\Lnorm{2}}, 
	\end{aligned} \\
	\label{dissipation-1}
	&
	\begin{aligned}
	&  D_1 = D_1(t) := \norm{\dt^2 v_{\varepsilon}, \dt^2(\varepsilon w_{\varepsilon})}{\Lnorm{2}} + \norm{\dt v_{\varepsilon}, \dt (\varepsilon w_{\varepsilon})}{\Hnorm{2}} \\
	& ~~~~ ~~~~ + \norm{\dt^2 \sigma_{\varepsilon}, \dt \nablah \sigma_{\varepsilon} , \dfrac{\dz \dt \sigma_{\varepsilon}}{\varepsilon}}{\Hnorm{1}} + \norm{\dt w_{\varepsilon}, \dz \dt w_{\varepsilon} }{\Hnorm{1}}.
	\end{aligned}
	\end{align}
	
	We will use $ \mathfrak N (\cdot)  $ to denote a locally Lipschitz nonlinear function of its argument(s), which can be different from line to line. For any two quantities $ Q_1 $ and $ Q_2 $, $ Q_1 \lesssim Q_2 $ is used to represent $ Q_1 \leq C Q_2 $ for some constant $ C \in (0,\infty) $, whose value will be different from line to line.
	
	\subsection{Main theorems}\label{subsec:theorem}
	We consider \eqref{CF} with initial data
	\begin{equation}\label{initial:CF}
			(\sigma_\varepsilon, v_\varepsilon, w_\varepsilon)\big|_{t=0} = ( \sigma_{0,\varepsilon},v_{0,\varepsilon}, w_{0,\varepsilon} ).
	\end{equation}
	Then the initial data for the time derivatives $ \dt \sigma_\varepsilon, \dt v_\varepsilon, \dt w_\varepsilon, \partial_t^2 \sigma_\varepsilon $ are given through compatibility by employing the equations in \eqref{CF}, inductively, i.e.,
	\begin{equation*}
		(\dt \sigma_\varepsilon, \dt v_\varepsilon, \dt w_\varepsilon, \partial_t^2 \sigma_\varepsilon)\big|_{t=0} = (\sigma_{1,\varepsilon},v_{1,\varepsilon},w_{1,\varepsilon},\sigma_{2,\varepsilon}),
	\end{equation*}
	where $ \sigma_{1,\varepsilon},v_{1,\varepsilon},w_{1,\varepsilon},\sigma_{2,\varepsilon},v_{2,\varepsilon},w_{2,\varepsilon} $ are given by
	\begin{gather*}
		\sigma_{1,\varepsilon} + v_{0,\varepsilon}\cdot\nablah \sigma_{0,\varepsilon} + w_{0,\varepsilon} \dz \sigma_{0,\varepsilon} + \dvh v_{0,\varepsilon} + \dz w_{0,\varepsilon} = 0,  \\
		v_{1,\varepsilon} + v_{0,\varepsilon} \cdot \nablah v_{0,\varepsilon} + w_{0,\varepsilon} \dz v_{0,\varepsilon} + \nablah \sigma_{0,\varepsilon} = \deltah v_{0,\varepsilon} + \partial_{zz} v_{0,\varepsilon},\\
		w_{1,\varepsilon} + v_{0,\varepsilon} \cdot \nablah w_{0,\varepsilon} + w_{0,\varepsilon} \dz w_{0,\varepsilon}  + \dfrac{1}{\varepsilon^2} \dz \sigma_{0,\varepsilon} =  \deltah w_{0,\varepsilon} + \partial_{zz} w_{0,\varepsilon}, \\
		\sigma_{2,\varepsilon} + v_{0,\varepsilon}\cdot\nablah \sigma_{1,\varepsilon} + w_{0,\varepsilon} \dz \sigma_{1,\varepsilon}+ v_{1,\varepsilon}\cdot\nablah \sigma_{0,\varepsilon} \\
		 + w_{1,\varepsilon} \dz \sigma_{0,\varepsilon} + \dvh v_{1,\varepsilon} + \dz w_{1,\varepsilon} = 0.
	\end{gather*}

	The first theorem in this paper is concerning the uniform-in-$\varepsilon$ estimates and the justification of hydrostatic approximation limit:
	\begin{thm}[Hydrostatic approximation]\label{thm:unifrom_hydrostatic}
		Suppose that the initial data in \eqref{initial:CF} satisfies $ E(0) < \infty $. Then there is a $ T^* \in (0,\infty) $, independent of $ \varepsilon $, such that the solution to \eqref{CF} with initial data \eqref{initial:CF} satisfies
		\begin{equation}\label{thm-001}
			\sup_{0\leq t \leq T^* } E^2(t) + \int_0^{T^*} D^2(t)\,dt < C_0 < \infty, 
		\end{equation}
		where $ C_0 $ is some positive constant depending only on $ E(0) $ and is independent of $ \varepsilon $. Moreover, there exist $ (\sigma_p, v_p, w_p )  $ with
\begin{equation*}
\begin{gathered}
\sigma_p \in L^\infty(0,T^{*};H^4), ~~ \dt \sigma_p \in L^\infty(0,T^*;H^2) \cap L^2(0,T^*;H^3), \\
v_p\in L^\infty(0,T^*;H^3) \cap L^2(0,T^*;H^4), ~~ \dt v_p \in L^2(0,T^*;H^2), \\
w_p, \dz w_p \in L^\infty(0,T^*;H^2)\cap L^2(0,T^*;H^3),
\end{gathered}
\end{equation*}
such that for a subsequence of $\lbrace (\sigma_\varepsilon, v_\varepsilon, w_\varepsilon) \rbrace$, as $ \varepsilon \rightarrow 0^+ $, 
\begin{align*}
 \sigma_\varepsilon & \buildrel\ast\over\rightharpoonup \sigma_p & ~ \text{weak-$\ast$ in} ~ &  L^\infty(0,T^*;H^4), \\
 \sigma_\varepsilon & \rightarrow \sigma_p & ~ \text{in} ~  & L^\infty(0,T^*;H^3) \cap C([0,T^*];H^3),\\
 \dt \sigma_\varepsilon, w_\varepsilon, \dz w_\varepsilon & \buildrel\ast\over\rightharpoonup \dt \sigma_p, w_p, \dz w_p & ~ \text{weak-$\ast$ in} ~ & L^\infty(0,T^*;H^2), \\
 \dt \sigma_\varepsilon , w_\varepsilon, \dz w_\varepsilon & \rightharpoonup \dt \sigma_p, w_p, \dz w_p & ~ \text{weakly in} ~ & L^2(0,T^*;H^3), \\
 v_\varepsilon & \buildrel\ast\over\rightharpoonup v_p & ~ \text{weak-$\ast$ in} ~ & L^\infty(0,T^*;H^3), \\
  v_\varepsilon & \rightarrow v_p & ~ \text{in} ~ & L^\infty(0,T^*;H^2)\cap C([0,T^*];H^2), \\
   v_\varepsilon & \rightharpoonup v_p & ~ \text{weakly in} ~ & L^2(0,T^*;H^4), \\
     \dt v_\varepsilon & \rightharpoonup \dt v_p & ~ \text{weakly in} ~ & L^2(0,T^*;H^2),
	\end{align*}
	and $ (\sigma_p, v_p, w_p) $ is a solution to \eqref{CPE}.

	In addition, suppose that the initial data in \eqref{initial:CF} satisfies $ E_1(0) < \infty $. Then there is a $ T^{**} \in (0,T^*] $, independent of $ \varepsilon $, such that 
		\begin{equation}\label{thm-002}
			\sup_{0\leq t \leq T^{**} } ( E^2(t) + E_1^2(t) ) + \int_0^{T^{**}} ( D^2(t) + D_1(t)) \,dt < C_1 < \infty, 
			\end{equation}
		where $ C_1 $ is some positive constant depending only on $ E(0) $ and $ E_1(0) $ and is independent of $ \varepsilon $. 
	\end{thm}

	We summarize the convergence rates in the following theorem:
	\begin{thm}[Rates of convergence]\label{thm:rate}
		Suppose that the solution $ (\sigma_p, v_p, w_p) $ to \eqref{CPE} given by the limit in Theorem \ref{thm:unifrom_hydrostatic} satisfies
		\begin{gather*}
			\norm{\sigma_p}{L^\infty(0,T^{**};\Hnorm{4})} + \norm{\dt \sigma_p}{L^\infty(0,T^{**};\Hnorm{2})}+ \norm{v_p}{L^\infty(0,T^{**};\Hnorm{3})}  \\
			 + \norm{v_p}{L^2(0,T^{**};\Hnorm{4})} + \norm{w_p}{L^\infty(0,T^{**};\Hnorm{2})} < C_p < \infty,
		\end{gather*}
		for some constant $ C_p \in (0,\infty) $,
		and the initial data in \eqref{initial:CF} satisfy
		\begin{equation}\label{WP}
			\norm{\sigma_{0,\varepsilon} - \sigma_p\big|_{t=0} , v_{0,\varepsilon} - v_p\big|_{t=0}}{\Lnorm{2}} \lesssim \varepsilon. 
		\end{equation}
		Then under the conditions in Theorem \ref{thm:unifrom_hydrostatic}, i.e., $ E(0), E_1(0) < \infty $, we have
		\begin{equation}\label{thm-003}
			\begin{gathered}
			\norm{ \sigma_\varepsilon - \sigma_p, v_\varepsilon - v_p }{L^\infty(0,T^{**};\Lnorm{2})} + \norm{v_\varepsilon - v_p}{L^2(0,T^{**};\Hnorm{1})} \leq C_2 \varepsilon, \\
			\norm{w_\varepsilon - w_p}{\Lnorm{\infty}(0,T^{**};L^2)} \leq C_2 \varepsilon^{2/3}, ~~ \norm{w_\varepsilon - w_p}{\Lnorm{2}(0,T^{**};L^2)} \leq C_2 \varepsilon^{3/4},
			\end{gathered}
		\end{equation}
		where $ C_2 \in (0,\infty) $ is a constant depending only on $ E(0), E_1(0) $ and $ C_p $, and is independent of $ \varepsilon $.
	\end{thm}
	
	\begin{proof}[Proof of Theorem \ref{thm:unifrom_hydrostatic}]
		The uniform estimates in \eqref{thm-001} and \eqref{thm-002} are  shown in section \ref{subsec:uniform-est} and section \ref{subsec:uniform-est-2}, given by \eqref{uniform-est-2} and \eqref{uniform-est-3}, respectively. The convergence is given in section \ref{sec:convergence-cpe}, below.
	\end{proof}
	
	\begin{proof}[Proof of Theorem \ref{thm:rate}]
		\eqref{thm-003} is the consequence of \eqref{con-rate-dlsigma-dlv-3} and \eqref{con-rate-204}, in section \ref{subsec:converging-rate}, below. 
	\end{proof}

	\begin{rmk}
		Assumption \eqref{WP} on the initial data above is the definition of well-prepared initial data, i.e., it essentially assumes that the initial data is close to the hydrostatic approximation. 
	\end{rmk}
	
	\begin{rmk}
		Theorem \ref{thm:unifrom_hydrostatic} guarantees the convergence of a subsequence to a solution of the limit equations. However, since the strong solution to the limit equations is well-posed, and in particular, is unique (see, e.g., \cite{LT2018a}), the convergence is actually of the full sequence. 
		Theorem \ref{thm:rate} 
		states the convergence rate of the full sequence for well-prepared initial data.
	\end{rmk}

\section{Uniform-in-$\varepsilon$ estimates}\label{sec:uniform-est}
	In this section, we shorten the notations by dropping the subscript $ \varepsilon $ in $ (\sigma_\varepsilon, v_\varepsilon, w_\varepsilon) $, i.e., $ (\sigma, v, w) = (\sigma_\varepsilon, v_\varepsilon, w_\varepsilon) $. 
Recall that we already have short time existence and uniqueness of solutions to \eqref{CF} on a time interval that might depend on $ \varepsilon $ (see, e.g., \cite{Cho2004}). The main goal of this section is to obtain the existence time and corresponding estimates that are independent of $ \varepsilon $. 
	
	\subsection{$ L^2 $--estimates}
	Take the $ L^2 $-inner product of \subeqref{CF}{2} with $ \dt v - \Delta v $, and \subeqref{CF}{3} with $  \dt w - \Delta w  $, respectively. One obtains, 
	\begin{align}
	& \label{est-001}
	\begin{aligned}
	& \norm{\dt v - \Delta v }{\Lnorm{2}}^2 = - \int \bigl( \dt v - \Delta v \bigr) 
	\cdot \bigl( \nablah \sigma + v \cdot \nablah v + w \dz v   \bigr)  \idx \\
	& ~~ \leq \dfrac{1}{2} \norm{\dt v - \Delta v}{\Lnorm{2}}^2 + C \bigl( \norm{\nablah \sigma}{\Lnorm{2}}^2 + \norm{v \cdot \nablah v}{\Lnorm{2}}^2 + \norm{w \dz v}{\Lnorm{2}}^2 \bigr),
	\end{aligned}\\
	& \label{est-002}
	\begin{aligned}
	& \varepsilon^2 \norm{\dt w - \Delta w }{\Lnorm{2}}^2 = - \varepsilon^2 \int \bigl( \dt w - \Delta w \bigr)
	\bigl( \dfrac{1}{\varepsilon^2} \dz \sigma + v\cdot \nablah w + w \dz w  \bigr) \idx \\
	& ~~ \leq \dfrac{\varepsilon^2}{2}  \norm{\dt w - \Delta w}{\Lnorm{2}}^2 
	+ \dfrac{C}{\varepsilon^2} \norm{\dz \sigma}{\Lnorm{2}}^2 + C \varepsilon^2 \bigl( \norm{v \cdot \nablah w}{\Lnorm{2}}^2 
	 + \norm{ w \dz w}{\Lnorm{2}}^2 \bigr),
	\end{aligned}
	\end{align}
	where we have applied the H\"older and Young inequalities. 
	Notice that, for $ \eta = v, w $, applying integration by parts yields
	\begin{equation}\label{id:int-by-parts-001}
	\norm{\dt \eta - \Delta \eta }{\Lnorm{2}}^2 = \dfrac{d}{dt} \norm{\nabla \eta }{\Lnorm{2}}^2 + \norm{\dt \eta}{\Lnorm{2}}^2 + \norm{\nabla^2 \eta}{\Lnorm{2}}^2.
	\end{equation}
	Thus \eqref{est-001} and \eqref{est-002} can be written as, 
	\begin{equation}\label{est-003}
	\begin{aligned}
	& \dfrac{d}{dt}  \norm{\nabla v,\varepsilon \nabla w }{\Lnorm{2}}^2  + \norm{ \dt v, \nabla^2 v,\varepsilon\dt w,  \varepsilon \nabla^2 w }{\Lnorm{2}}^2 \leq C \norm{\nablah \sigma, \dfrac{\dz \sigma}{\varepsilon}}{\Lnorm{2}}^2 \\
	& ~~~~ ~~~~ + C \sum_{i=1}^2 \norm{\mathcal I_i}{\Lnorm{2}}^2, 
	\end{aligned}\end{equation}
	where
	\begin{equation}\label{id:nonlinearity-1}
	\begin{aligned}
	& \mathcal I_1 = \bigl( v \cdot\nablah v, v \cdot \nablah (\varepsilon w) \bigr), ~~~
	\mathcal I_2 =  \bigl( w\dz v, w \dz (\varepsilon w) \bigr). 
	\end{aligned}
	\end{equation}
	
	We postpone the estimates for $ \norm{\mathcal I_i}{\Lnorm{2}} $ to section \ref{subsec:est-nonlinearities}. However, even without detailed calculations, one can see that, a uniform-in-$\varepsilon$ estimate for $ w $ is necessary, for instance, to get an estimate for $ \mathcal I_2 $. This indeed is the main challenge to be addressed in this work. In the next section, we discuss our strategy to overcome this obstacle.

	\subsection{Estimates for the vertical velocity and the density}
	
	We use \subeqref{CF}{1} to represent $ w $ in terms of $ \sigma $ and $ v $. Indeed, after multiplying \subeqref{CF}{1} with $ e^{\sigma} $, it follows that
	\begin{equation*}
	\dt e^\sigma + \dvh (e^{\sigma} v) + \dz (e^{\sigma} w) = 0.
	\end{equation*}
	Thus, recalling that $ \cdot_h $ represents the horizontal variables $x,y $, one has that,
	\begin{equation}\label{id:vertical-velocity-0}
\begin{gathered}
	w(\cdot_h,z) = - e^{-\sigma(\cdot_h,z)} \int_0^z e^{\sigma(\cdot_h,z')} \biggl( \dt \sigma(\cdot_h,z') + v(\cdot_h,z') \cdot \nablah \sigma(\cdot_h,z') \\ + \dvh v(\cdot_h,z') \biggr) \,dz'.
\end{gathered}
	\end{equation}
	Similarly, we remark that we can represent $ \dlw $ in terms of $ \dlsigma $, $ \dlv $, $ \sigma_p $, $ v_p $, and $ w_p $, by using \subeqref{eq:difference}{1}, even thought we don't need it now. One can obtain, after multiplying \subeqref{eq:difference}{1} with $ e^{\dlsigma} $, that
	\begin{equation*}
	\dt e^{\dlsigma} + v_p \cdot \nablah e^{\dlsigma} + w_p \dz e^{\dlsigma} + e^{\dlsigma}  \dlv \cdot \nablah \sigma_p + \dvh (e^{\dlsigma} \dlv) + \dz (e^{\dlsigma}\dlw) = 0,
	\end{equation*}
	and
	\begin{equation}\label{id:vertical-velocity-1}
	\begin{aligned}
	& \dlw(\cdot_h,z) =  - e^{-\dlsigma(\cdot_h,z)} \int_0^z e^{\dlsigma(\cdot_h,z')} \biggl(\dt \dlsigma(\cdot_h,z') + v(\cdot_h,z') \cdot \nablah \dlsigma(\cdot_h,z')\\
	& ~~ + \dlv(\cdot_h,z') \cdot \nablah \sigma_p(\cdot_h,z') + w_p(\cdot_h,z') \dz \dlsigma (\cdot_h,z') 
	+ \dvh \dlv(\cdot_h,z') \biggr) \,dz'.
	\end{aligned}
	\end{equation}
Eventually, we will return to \eqref{id:vertical-velocity-1}. 
	Then directly using the representation of $ w $ in \eqref{id:vertical-velocity-0}, we have the following:
	
	\begin{prop}\label{lm:vertical-velocity}
	Assuming $ w $ is given by \eqref{id:vertical-velocity-0}, one can show that
		\begin{align}
		\label{est-w-1}
		& \norm{w, \dz w}{\Hnorm{2}} \leq C e^{2 \norm{\sigma}{\Hnorm{2}}} \bigl( \norm{\sigma}{\Hnorm{3}}^3 + 1 \bigr) \bigl( \norm{\dt\sigma}{\Hnorm{2}} + \norm{v}{\Hnorm{3}} + \norm{\sigma}{\Hnorm{3}} \norm{v}{\Hnorm{2}}  \bigr) ,  \\
		\label{est-w-2}
		& \norm{w, \dz w}{\Hnorm{3}} \leq C e^{2 \norm{\sigma}{\Hnorm{2}}} \bigl( \norm{\sigma}{\Hnorm{4}}^4 + 1 \bigr) \bigl( \norm{\dt\sigma}{\Hnorm{3}} + \norm{v}{\Hnorm{4}} + \norm{\sigma}{\Hnorm{4}} \norm{v}{\Hnorm{3}}  \bigr),
		\end{align}
		for some generic positive constant $ C \in (0,\infty) $, independent of $ \varepsilon $. 
	\end{prop}
	\begin{proof}
		To simplify the notations, denote by
		\begin{equation}\label{id:vertical-velocity-notation}
		\Xi : = \dt \sigma + v \cdot \nablah \sigma + \dvh v.
		\end{equation}
		Next, we calculate the derivatives of $ w $:
		\begin{align*}
		& \dz w = e^{-\sigma} \dz \sigma \int_0^z \bigl( e^{\sigma} \Xi \bigr)  \,dz' - \Xi, \\
		& \ddh w = e^{-\sigma} \ddh \sigma \int_0^z \bigl( e^{\sigma} \Xi \bigr) \,dz' 
		- e^{-\sigma}\int_0^z \bigl( e^{\sigma} \ddh \sigma \Xi  + e^{\sigma} \ddh \Xi\bigr) \,dz', \\
		& \partial_{zz} w = e^{-\sigma}\bigl(\partial_{zz}\sigma - (\dz \sigma)^2 \bigr)\int_0^z \bigl( e^{\sigma} \Xi  \bigr) \,dz' 
		+ \dz \sigma \Xi  - \dz \Xi , \\
		& \partial_{hz} w = e^{-\sigma}\bigl( \partial_{hz} \sigma - \ddh\sigma \dz \sigma \bigr) \int_0^z \bigl( e^{\sigma} \Xi \bigr) \,dz' \\
		& ~~~~ + e^{-\sigma}\dz\sigma \int_0^z \bigl( e^{\sigma} \ddh \sigma \Xi  + e^{\sigma} \ddh \Xi \bigr) \,dz' 
		- \ddh \Xi,\\
		& \partial_{hh} w = e^{-\sigma} \bigl(\partial_{hh}\sigma - (\ddh \sigma )^2 \bigr)\int_0^z \bigl( e^{\sigma} \Xi  \bigr) \,dz' \\
		& ~~~~ + 2 e^{-\sigma}\ddh \sigma \int_0^z \bigl( e^{\sigma} \ddh \sigma \Xi  
		+ e^{\sigma} \ddh \Xi  \bigr) \,dz'\\
		& ~~~~ - e^{-\sigma} \int_0^z \bigl\lbrack e^{\sigma}\bigl(\partial_{hh}\sigma + (\ddh \sigma)^2 \bigr) \Xi 
		+ 2 e^{\sigma} \ddh \sigma \ddh \Xi 
		+ e^\sigma \partial_{hh}\Xi  \bigr\rbrack \,dz',\\
		& \partial_{zzz} w = e^{-\sigma} \bigl( \partial_{zzz}\sigma - 3 \partial_z\sigma \partial_{zz}\sigma + (\dz\sigma)^3 \bigr) \int_0^z \bigl( e^{\sigma}\Xi \bigr) \,dz' \\
		& ~~~~ + \bigl(2 \partial_{zz} \sigma - (\partial_z\sigma)^2 \bigr) \Xi
		+ \dz \sigma \dz \Xi  - \partial_{zz}\Xi , \\
		& \partial_{hzz} w =  e^{-\sigma} \bigl( \partial_{hzz} \sigma - 2 \partial_{hz}\sigma \dz \sigma - \ddh \sigma \partial_{zz}\sigma  + \ddh \sigma (\dz\sigma)^2 \bigr) \\
		& ~~~~ ~~~~ \times  \int_0^z \bigl( e^{\sigma} \Xi \bigr) \,dz' + e^{-\sigma} \bigl( \partial_{zz} \sigma - (\dz \sigma)^2 \bigr)  \int_0^z \bigl( e^{\sigma}\ddh \sigma \Xi + e^{\sigma} \ddh \Xi \bigr) \,dz'  \\
		& ~~~~ + \partial_{hz} \sigma \Xi  
		+ \dz \sigma \ddh \Xi  - \partial_{hz} \Xi , \\
		& \partial_{hhz} w = e^{-\sigma} \bigl(\partial_{hhz}\sigma - 2\ddh \sigma \partial_{hz}\sigma - \partial_{hh}\sigma \dz\sigma + (\ddh \sigma)^2 \dz \sigma  \bigr)\int_0^z \bigl( e^{\sigma} \Xi \bigr) \,dz' \\
		& ~~~~ + 2 e^{-\sigma} \bigl( \partial_{hz} \sigma - \ddh \sigma \dz \sigma \bigr) \int_0^z \bigl( e^{\sigma} \ddh \sigma \Xi  + e^{\sigma} \ddh \Xi \bigr) \,dz'\\
		& ~~~~ + e^{-\sigma} \dz \sigma \int_0^z \bigl\lbrack e^{\sigma}\bigl(\partial_{hh}\sigma + (\ddh \sigma)^2 \bigr)  \Xi 
		+ 2 e^{\sigma} \ddh \sigma  \ddh \Xi  
		+ e^\sigma \partial_{hh}\Xi \bigr\rbrack \,dz' \\
		& ~~~~ - \partial_{hh}\Xi,\\
		& \partial_{hhh} w = e^{-\sigma} \bigl( \partial_{hhh}\sigma - 3 \ddh \sigma \partial_{hh}\sigma + (\ddh \sigma)^3 \bigr) \int_0^z \bigl(e^\sigma \Xi\bigr)\,dz' \\
		& ~~~~ + 3 e^{-\sigma} \bigl( \partial_{hh}\sigma - (\ddh \sigma)^2 \bigr) \int_0^z \bigl( e^\sigma \ddh \sigma \Xi + e^\sigma \ddh \Xi \bigr) \,dz' \\
		& ~~~~ + 3 e^{-\sigma} \ddh \sigma \int_0^z \bigl\lbrack e^\sigma \bigl(\partial_{hh}\sigma + (\ddh \sigma)^2 \bigr) \Xi + 2 e^\sigma \ddh \sigma \ddh \Xi + e^\sigma \partial_{hh}\Xi \bigr\rbrack\,dz'\\
		& ~~~~ - e^{-\sigma} \int_0^z \bigl\lbrack e^\sigma\bigl( \partial_{hhh}\sigma + 3 \ddh \sigma \partial_{hh}\sigma + (\ddh \sigma)^3 \bigr)\Xi  + 3 e^\sigma \bigl(\partial_{hh}\sigma + (\ddh \sigma)^2 \bigr) \ddh \Xi\\
		& ~~~~ ~~~~ + 3 e^\sigma \ddh \sigma \partial_{hh}\Xi + e^\sigma \partial_{hhh} \Xi \bigr\rbrack 
		\,dz',\\
		& \partial_{zzzz} w = e^{-\sigma} \bigl\lbrack \partial_{zzzz} \sigma - 3 (\partial_{zz} \sigma)^2 - 4 \dz \sigma \partial_{zzz}\sigma + 6 (\dz \sigma
		)^2 \partial_{zz}\sigma - (\dz \sigma)^4 \bigr\rbrack \\
		& ~~~~ ~~~~ \times \int_0^z \bigl( e^\sigma \Xi \bigr) \,dz' + \bigl( 3\partial_{zzz}\sigma - 5 \dz \sigma \partial_{zz}\sigma + (\dz \sigma)^3 \bigr) \Xi \\
		& ~~~~ + \bigl( 3 \partial_{zz}\sigma - (\dz \sigma)^2\bigr)\dz \Xi 
		+ \dz \sigma \partial_{zz}\Xi - \partial_{zzz}\Xi, \\
		& \partial_{hzzz} w = e^{-\sigma} \bigl\lbrack \partial_{hzzz} \sigma - 3 \partial_{hzz}\sigma \dz\sigma - 3 \partial_{hz}\sigma \partial_{zz}\sigma - \ddh \sigma \partial_{zzz} \sigma \\
		& ~~~~ ~~~~ + 3 \partial_{hz}\sigma (\dz \sigma)^2 + 3 \ddh \sigma \dz\sigma \partial_{zz}\sigma - \ddh \sigma(\dz \sigma)^3 \bigr\rbrack \times \int_0^z \bigl( e^\sigma \Xi \bigr) \,dz'  \\
		& ~~~~ + e^{-\sigma}\bigl\lbrack \partial_{zzz}\sigma - 3 \dz \sigma \partial_{zz} \sigma + (\dz \sigma)^3 \bigr\rbrack \int_0^z\bigl(e^\sigma \ddh \sigma \Xi + e^\sigma \ddh \Xi \bigr) \,dz' \\
		& ~~~~ + \bigl( 2 \partial_{hzz}\sigma - 2 \partial_{hz}\sigma\dz \sigma \bigr) \Xi + \bigl( 2 \partial_{zz}\sigma - (\dz \sigma)^2 \bigr) \ddh \Xi + \partial_{hz}\sigma \dz \Xi\\
		& ~~~~ + \dz \sigma \partial_{hz}\Xi - \partial_{hzz}\Xi,  \\
		& \partial_{hhzz} w = e^{-\sigma} \bigl\lbrack \partial_{hhzz} \sigma - 2 (\partial_{hz}\sigma)^2 - 2 \ddh \sigma \partial_{hzz}\sigma - 2 \partial_{hhz}\sigma \dz \sigma \\
		& ~~~~ ~~~~ - \partial_{hh}\sigma \partial_{zz}\sigma  + 4 \ddh \sigma \dz \sigma \partial_{hz}\sigma + (\ddh \sigma)^2 \partial_{zz}\sigma + \partial_{hh}\sigma (\dz \sigma)^2  \\
		& ~~~~ ~~~~ - (\ddh \sigma)^2 (\dz \sigma)^2 \bigr\rbrack
		\times \int_0^z \bigl( e^\sigma \Xi \bigr) \,dz'
		+ 2 e^{-\sigma} \bigl\lbrack \partial_{hzz}\sigma - 2 \partial_{hz}\sigma \dz\sigma \\
		& ~~~~ ~~~~ - \ddh \sigma \partial_{zz}\sigma + \ddh \sigma(\dz \sigma)^2 \bigr\rbrack \times \int_0^z \bigl(e^\sigma \ddh \sigma \Xi + e^\sigma \ddh \Xi \bigr) \,dz' \\
		& ~~~~ + e^{-\sigma} \bigl\lbrack \partial_{zz}\sigma - (\dz \sigma)^2 \bigr\rbrack \int_0^z \bigl\lbrack e^\sigma \bigl(\partial_{hh}\sigma + (\ddh \sigma)^2 \bigr) \Xi + 2 e^\sigma \ddh \sigma \ddh \Xi \\
		& ~~~~ ~~~~ + e^\sigma \partial_{hh}\Xi \bigr\rbrack\,dz' 
		+ \partial_{hhz}\sigma \Xi + 2 \partial_{hz}\sigma \ddh \Xi + \dz\sigma \partial_{hh}\Xi - \partial_{hhz}\Xi, \\
		& \partial_{hhhz} w = e^{-\sigma} \bigl\lbrack \partial_{hhhz} \sigma - 3 \partial_{hz} \sigma \partial_{hh}\sigma - 3 \ddh \sigma \partial_{hhz}\sigma + 3 (\ddh \sigma)^2 \partial_{hz}\sigma \\
		& ~~~~ ~~~~ - \partial_{hhh} \sigma \dz \sigma + 3 \ddh \sigma \dz \sigma \partial_{hh}\sigma - (\ddh \sigma)^3 \dz \sigma  \bigr\rbrack \times \int_0^z \bigl(e^\sigma \Xi \bigr) \,dz' \\
		& ~~~~ + 3 e^{-\sigma} \bigl\lbrack \partial_{hhz}\sigma - 2 \ddh \sigma \partial_{hz}\sigma - \partial_{hh}\sigma \dz \sigma + (\ddh \sigma)^2 \dz \sigma \bigr\rbrack \\
		& ~~~~ ~~~~ \times \int_0^z \bigl( e^\sigma \ddh \sigma \Xi + e^\sigma \ddh \Xi \bigr)\,dz' + 3 e^{-\sigma} \bigl\lbrack \partial_{hz}\sigma - \ddh \sigma \dz \sigma \bigr\rbrack \\
		& ~~~~ ~~~~ \times \int_0^z \bigl\lbrack e^\sigma \bigl( \partial_{hh}\sigma+(\ddh\sigma)^2 \bigr) \Xi + 2 e^\sigma \ddh \sigma \ddh \Xi + e^\sigma \partial_{hh}\Xi  \bigr\rbrack\,dz'\\
		& ~~~~ + e^{-\sigma} \dz \sigma \int_0^z \bigl\lbrack e^\sigma\bigl( \partial_{hhh}\sigma +3 \ddh \sigma \partial_{hh}\sigma +(\ddh\sigma)^3 \bigr) \Xi  \\
		& ~~~~ ~~~~ + 3 e^\sigma \bigl( \partial_{hh}\sigma + (\ddh \sigma)^2 \bigr) \ddh \Xi + 3 e^\sigma\ddh \sigma \partial_{hh}\Xi + e^\sigma\partial_{hhh}\Xi  \bigr\rbrack \,dz' - \partial_{hhh}\Xi.
		\end{align*}

		Next, after applying the H\"older, Minkowski, the Sobolev embedding inequalities, and that $ H^m $ is an algebra for $ m \geq 2 $, we obtain
		\begin{equation}
		\label{est-w-001}
		\begin{aligned}
		& \norm{w}{\Lnorm{2}} \lesssim \norm{e^{-\sigma}}{\Lnorm{\infty}}\norm{e^{\sigma} \bigl( \dt \sigma + v \cdot \nablah \sigma + \dvh v \bigr)}{\Lnorm{2}} \\
		& ~~~~ \lesssim \norm{e^{-\sigma}}{\Lnorm{\infty}} \norm{e^{\sigma}}{\Lnorm{\infty}} \bigl( \norm{\dt \sigma}{\Lnorm{2}} + \norm{v}{\Lnorm{6}}\norm{\nablah \sigma}{\Lnorm{3}} + \norm{\nablah v}{\Lnorm{2}})\\
		& ~~~~ \lesssim e^{2\norm{\sigma}{\Hnorm{2}}} \bigl(\norm{\dt \sigma}{\Lnorm{2}} + \norm{v}{\Hnorm{1}}  \norm{\sigma}{\Hnorm{2}} + \norm{v}{\Hnorm{1}} \bigr).
		\end{aligned}
		\end{equation}
		Similarly, with details omitted, one can also show,
		\begin{align}
		& \label{est-w-002}
		\begin{aligned}
		& \norm{\dz w, \nablah w, \partial_{zz} w, \nablah \partial_{z} w}{\Lnorm{2}} \lesssim e^{2\norm{\sigma}{\Hnorm{2}}}\bigl( \norm{\sigma}{\Hnorm{3}}^2 + 1\bigr) \\
		& ~~~~ ~~~~\times \bigl( \norm{\dt\sigma}{\Hnorm{1}} + \norm{\sigma}{\Hnorm{3}}\norm{v}{\Hnorm{1}} + \norm{v}{\Hnorm{2}} \bigr),
		\end{aligned}\\
		& \label{est-w-003}
		\begin{aligned}
		& \norm{\nablah^2 w}{\Lnorm{2}} \lesssim e^{2\norm{\sigma}{\Hnorm{2}}}\bigl( \norm{\sigma}{\Hnorm{3}}^2 + 1\bigr) \bigl( \norm{\dt\sigma}{\Hnorm{2}} + \norm{\sigma}{\Hnorm{3}}\norm{v}{\Hnorm{2}} + \norm{v}{\Hnorm{3}} \bigr),
		\end{aligned}\\
		& \label{est-w-004}
		\begin{aligned}
		& \norm{\nabla^2 \dz w}{\Lnorm{2}} \lesssim e^{2\norm{\sigma}{\Hnorm{2}}}\bigl( \norm{\sigma}{\Hnorm{3}}^3 + 1\bigr) \bigl( \norm{\dt\sigma}{\Hnorm{2}} + \norm{\sigma}{\Hnorm{3}}\norm{v}{\Hnorm{2}} + \norm{v}{\Hnorm{3}} \bigr),
		\end{aligned}\\
		& \label{est-w-005}
		\begin{aligned}
		& \norm{\nablah^3 w}{\Lnorm{2}}\lesssim e^{2\norm{\sigma}{\Hnorm{2}}}\bigl( \norm{\sigma}{\Hnorm{3}}^3 + 1\bigr) \bigl( \norm{\dt\sigma}{\Hnorm{3}} + \norm{\sigma}{\Hnorm{4}}\norm{v}{\Hnorm{3}} + \norm{v}{\Hnorm{4}} \bigr), 
		\end{aligned} \\
		& \label{est-w-006}
		\begin{aligned}
		& \norm{\nabla^3 \dz w}{\Lnorm{2}}\lesssim e^{2\norm{\sigma}{\Hnorm{2}}}\bigl( \norm{\sigma}{\Hnorm{4}}^4 + 1\bigr) \bigl( \norm{\dt\sigma}{\Hnorm{3}} + \norm{\sigma}{\Hnorm{4}}\norm{v}{\Hnorm{3}} + \norm{v}{\Hnorm{4}} \bigr).
		\end{aligned}
		\end{align}
		We summarize inequalities \eqref{est-w-1} and \eqref{est-w-2} from inequalities \eqref{est-w-001}--\eqref{est-w-006}.
	\end{proof}

	In view of \eqref{est-w-1} and \eqref{est-w-2}, the estimates for $ w $ can be bounded by the estimates for $ \dt \sigma $, which may not be bounded and induce fast oscillation as $ \varepsilon \rightarrow 0 ^+ $. Therefore, we will need to obtain the uniform estimates for $ \dt \sigma $. To achieve this goal, we formulate a strongly damped wave equation for $ \sigma $ in the following (see \eqref{eq:mixed}, below). In fact, after applying $ \dt $ to \subeqref{CF}{1}, $ \dvh $ to \subeqref{CF}{2}, and $ \dz $ to \subeqref{CF}{3}, and combining the resultant equations, we obtain, with direct calculations,
	\begin{equation}\label{eq:pre-mixed-1}
	\begin{aligned}
	& \partial_{tt} \sigma + v \cdot \nablah \dt \sigma + w \dz \dt \sigma + \dt v \cdot \nablah \sigma + \dt w \dz \sigma \\
	& ~~ = - \dvh \dt v - \dz \dt w = \deltah \sigma + \dfrac{1}{\varepsilon^2} \partial_{zz} \sigma - \Delta ( \dvh v + \partial_z w) \\
	& ~~~~ + \dvh \bigl( v \cdot \nablah v + w \dz v \bigr) 
	+ \partial_z \bigl( v \cdot \nablah w + w \dz w \bigr).
	\end{aligned}
	\end{equation}
	On the other hand, from \subeqref{CF}{1}, we have
	\begin{align*}
	& - \Delta ( \dvh v + \partial_z w) = \Delta ( \dt \sigma + v\cdot \nablah \sigma 
	 + w\dz \sigma ) \\
	& = \dt \Delta \sigma 
	 + (v \cdot\nablah + w \dz) \Delta \sigma + \Delta v \cdot \nablah \sigma + 2 \nabla v : \nablah \nabla \sigma \\
	& ~~~~ + \Delta w \dz \sigma + 2 \nabla w \cdot \dz \nabla \sigma.
	\end{align*}
	Hence, after substituting the above identity to \eqref{eq:pre-mixed-1}, we arrive at the following equation:
	\begin{equation}\label{eq:mixed}
	\dt \bigl(\dt \sigma - \Delta \sigma) + (v \cdot\nablah + w \dz) \bigl( \dt \sigma - \Delta \sigma \bigr) - \deltah \sigma - \dfrac{1}{\varepsilon^2} \partial_{zz} \sigma 
	=  \sum_{j=1}^{4} \mathcal{J}_j,
	\end{equation}
	where
	\begin{equation}\label{id:nonlinearity-2}
	\begin{aligned}
	& \mathcal J_1 : = \dt v \cdot \nablah \sigma - \Delta v \cdot\nablah \sigma - 2 \nabla v : \nablah \nabla \sigma,\\
	& \mathcal J_2 := \dt w \dz \sigma - \Delta w \dz \sigma 
	- 2 \nabla w \cdot \dz \nabla \sigma,\\
	& \mathcal J_3 := \dvh \bigl( v \cdot \nablah v + w \dz v \bigr),\\
	& \mathcal J_4 : = \partial_z \bigl( v \cdot \nablah w + w \dz w  \bigr).
	\end{aligned}\end{equation}
	
	Next, we take the $ L^2 $--inner product of \eqref{eq:mixed} with $ \dt \sigma - \Delta\sigma $. Noticing that,
	\begin{align*}
	& \int \bigl( - \deltah \sigma - \dfrac{1}{\varepsilon^2}\partial_{zz} \sigma \bigr)\bigl(\dt \sigma - \Delta \sigma \bigr) \idx \\
	& ~~~~ = \dfrac{1}{2} \dfrac{d}{dt} \norm{\nablah \sigma, \dfrac{\dz \sigma}{\varepsilon}}{\Lnorm{2}}^2 + \norm{\nabla \nablah \sigma , \dfrac{\nabla \dz \sigma}{\varepsilon}}{\Lnorm{2}}^2,
	\end{align*}
	we arrive at:
	\begin{equation}\label{est-004}
	\begin{aligned}
	& \dfrac{1}{2} \dfrac{d}{dt} \norm{\dt \sigma - \Delta \sigma, \nablah \sigma, \dfrac{\dz \sigma}{\varepsilon}}{\Lnorm{2}}^2 + \norm{\nabla \nablah \sigma , \dfrac{\nabla \dz \sigma}{\varepsilon}}{\Lnorm{2}}^2 \\
	& ~~~~ = \dfrac{1}{2} \int \bigl(\dvh v + \partial_z w \bigr) |\dt \sigma - \Delta \sigma|^{2} \idx 
	+ \int \bigl(\dt \sigma - \Delta \sigma \bigr) \bigl(  \sum_{j=1}^{4} \mathcal{J}_j \bigr) \idx\\
	& ~~~~ \lesssim  \norm{\dt \sigma -\Delta \sigma}{\Lnorm{2}}^2\norm{\nablah v, \dz w}{\Hnorm{2}} + \norm{\dt \sigma -\Delta \sigma}{\Lnorm{2}} \sum_{j=1}^4 \norm{\mathcal J_j}{\Lnorm{2}}.
	\end{aligned}
	\end{equation}
	On the other hand, we take the $ L^2 $-inner product of \eqref{eq:mixed} with $ \dt \sigma $, and arrive at:
	\begin{equation}\label{est-005}
	\begin{aligned}
	& \dfrac{1}{2} \dfrac{d}{dt} \norm{\dt\sigma, \nablah \sigma, \dfrac{\dz \sigma}{\varepsilon}}{\Lnorm{2}}^2 + \norm{\nabla\dt \sigma}{\Lnorm{2}}^2 \\
	& ~~~~ = \dfrac{1}{2} \int (\dvh v + \partial_z w \bigr) \bigl( |\dt \sigma|^{2} 
	-2 \Delta \sigma \dt \sigma \bigr) \idx \\
	& ~~~~ - \int ( v \cdot \nablah + w \dz ) \dt \sigma \Delta \sigma \idx 
	+ \int \dt \sigma \bigl(  \sum_{j=1}^{4} \mathcal{J}_j\bigr) \idx \\
	& ~~~~ \lesssim \bigl( \norm{\dt \sigma}{\Lnorm{2}}^2 + \norm{\Delta \sigma}{\Lnorm{2}}^2 \bigr) \norm{\nablah v, \dz w}{\Hnorm{2}} \\
	& ~~~~ + \norm{v, w}{\Hnorm{2}} \norm{\nabla \dt \sigma}{\Lnorm{2}} \norm{\Delta \sigma}{\Lnorm{2}} + \norm{\dt \sigma}{\Lnorm{2}} \sum_{j=1}^4 \norm{\mathcal{J}_j}{\Lnorm{2}},
	\end{aligned}
	\end{equation}
	where we have substituted, after applying integration by parts,
	\begin{gather*}
	\int \biggl( (v\cdot \nablah + w \dz)\Delta \sigma \biggr) \dt \sigma \idx = - \int (\dvh v + \dz w) \Delta \sigma \dt \sigma \idx \\
	- \int \biggl( ( v \cdot \nablah + w \dz ) \dt \sigma \biggr) \Delta \sigma \idx.
	\end{gather*}
	Thus, \eqref{est-004} and \eqref{est-005} yield
	\begin{equation}\label{est-density-001}
	\begin{aligned}
	& \dfrac{d}{dt} \bigl\lbrace\norm{ \dt \sigma, \dt\sigma - \Delta \sigma}{\Lnorm{2}}^2 + 2 \norm{ \nablah \sigma, \dfrac{\dz \sigma}{\varepsilon}}{\Lnorm{2}}^2 \bigr\rbrace + 2 \norm{ \nabla \nablah \sigma, \dfrac{\nabla \dz \sigma}{\varepsilon}, \nabla \dt \sigma}{\Lnorm{2}}^2 \\
	&  \leq C \bigl( \norm{\dt\sigma}{\Lnorm{2}}^2 + \norm{\dt\sigma - \Delta \sigma}{\Lnorm{2}}^2) \norm{\nablah v, \dz w}{\Hnorm{2}}\\
	& ~~~~ + C \norm{v, w}{\Hnorm{2}} \norm{\nabla \dt \sigma}{\Lnorm{2}} \norm{\Delta \sigma}{\Lnorm{2}} \\
	& ~~~~
	 + C \bigl( \norm{\dt \sigma}{\Lnorm{2}} + \norm{\dt \sigma - \Delta \sigma}{\Lnorm{2}} \bigr) \sum_{j=1}^4 \norm{\mathcal{J}_j}{\Lnorm{2}}.
	\end{aligned}
	\end{equation}
	We postpone the estimates for $ \norm{\mathcal J_j}{\Lnorm{2}} $, $ j = 1,2,3,4 $, to section \ref{subsec:est-nonlinearities}. 

	\subsection{Higher order uniform-in-$\varepsilon$ estimates}
	In this section, we will derive higher order estimates, i.e., estimates for higher order spatial derivatives, corresponding to \eqref{est-003} and \eqref{est-density-001}. These are the essential parts for closing the estimates for the nonlinearities.

	Take the $ L^2 $-inner product of \subeqref{CF}{2} with $ \Delta^2 (\dt  v - \Delta v) $, and \subeqref{CF}{3} with $   \Delta^2 (\dt  w - \Delta w ) $, respectively. This yields, after applying integration by parts,
	\begin{align}
	& \label{est-201}
	\begin{aligned}
	& \norm{\Delta (\dt v - \Delta v) }{\Lnorm{2}}^2 = - \int \Delta (\dt v - \Delta v ) \cdot \Delta \bigl( \nablah \sigma + v \cdot \nablah v \\
	& ~~~~ ~~~~ + w \dz v  \bigr)\idx \leq \dfrac{1}{2} \norm{\Delta(\dt v - \Delta v)}{\Lnorm{2}}^2\\
	& ~~~~ + C \bigl( \norm{\Delta\nablah \sigma}{\Lnorm{2}}^2 + \norm{\Delta(v \cdot \nablah v)}{\Lnorm{2}}^2 + \norm{\Delta (w \dz v)}{\Lnorm{2}}^2 \bigr),
	\end{aligned} \\
	& \label{est-202}
	\begin{aligned}
	& \varepsilon^2 \norm{\Delta (\dt w - \Delta w) }{\Lnorm{2}}^2 = - \varepsilon^2 \int \Delta (\dt w - \Delta w) \times \Delta \bigl( \dfrac{1}{\varepsilon^2}  \dz \sigma \\
	& ~~~~ ~~~~ + v \cdot \nablah w + w \dz w \bigr) \idx 
	 \leq \dfrac{\varepsilon^2}{2} \norm{\Delta(\dt w - \Delta w)}{\Lnorm{2}}^2 \\
	& ~~~~
	+ \dfrac{C}{\varepsilon^2} \norm{\Delta \dz \sigma}{\Lnorm{2}}^2 + C \varepsilon^2 \bigl( \norm{\Delta (v\cdot \nablah w)}{\Lnorm{2}}^2 
	+ \norm{\Delta(w \dz w)}{\Lnorm{2}}^2 \bigr).
	\end{aligned}
	\end{align}
	Notice, similarly to \eqref{id:int-by-parts-001}, for $ \eta = v, w $, after applying integration by parts, we have the identity
	\begin{equation}\label{id:int-by-parts-002}
	\norm{\Delta(\dt \eta - \Delta \eta )}{\Lnorm{2}}^2 = \dfrac{d}{dt} \norm{\nabla^3 \eta}{\Lnorm{2}}^2 + \norm{\nabla^2\dt \eta}{\Lnorm{2}}^2 + \norm{\nabla^4 \eta}{\Lnorm{2}}^2.
	\end{equation}
	Therefore, \eqref{est-201} and \eqref{est-202} imply, after applying the H\"older inequality, that 
	\begin{equation}\label{est-203}
	\begin{aligned}
	& \dfrac{d}{dt} \norm{\nabla^3 v, \varepsilon \nabla^3 w}{\Lnorm{2}}^2 + \norm{\nabla^2 \dt v,  \nabla^4 v, \varepsilon \nabla^2 \dt w, \varepsilon \nabla^4 w}{\Lnorm{2}}^2  \\
	& ~~~~ \leq C \norm{\nabla^2 \nablah \sigma, \dfrac{\nabla^2 \dz \sigma}{\varepsilon}}{\Lnorm{2}}^2 + C \sum_{i=1}^{2} \norm{\mathcal I_i}{\Hnorm{2}}^2.
	\end{aligned}
	\end{equation}

	Notice, again, that after applying integration by parts, we have the following identities: 
	\begin{align*}
	& \int \dt (\dt \sigma - \Delta \sigma) \Delta^2 (\dt \sigma - \Delta \sigma) \idx = \dfrac{1}{2} \dfrac{d}{dt} \norm{\nabla^2 (\dt \sigma - \Delta \sigma) }{\Lnorm{2}}^2, \\
	& \int ( -\deltah \sigma - \dfrac{1}{\varepsilon^2}\partial_{zz}\sigma ) \Delta^2 (\dt \sigma - \Delta \sigma) \idx = \dfrac{1}{2} \dfrac{d}{dt} \norm{\nabla^2 \nablah \sigma, \dfrac{\nabla^2 \dz \sigma}{\varepsilon}  }{\Lnorm{2}}^2 \\
	& ~~~~ + \norm{ \nabla^3 \nablah \sigma, \dfrac{\nabla^3 \dz \sigma}{\varepsilon}}{\Lnorm{2}}^2, \\
	& \int (v \cdot \nablah + w \dz )\bigl(\dt \sigma - \Delta \sigma\bigr) \Delta^2 (\dt \sigma - \Delta \sigma)  \idx = \int  \Delta (\dt \sigma - \Delta \sigma) \\
	& ~~~~ ~~~~ \times \Delta\bigl\lbrack (v \cdot \nablah + w \dz )\bigl(\dt \sigma - \Delta \sigma\bigr) \bigr\rbrack  \idx \\
	& ~~~~ = \underbrace{\int  (v \cdot \nablah + w \dz ) \Delta (\dt \sigma - \Delta \sigma) \Delta (\dt \sigma - \Delta \sigma)  \idx}_{ = - \dfrac{1}{2} \int  (\dvh v + \dz w ) |\Delta (\dt \sigma - \Delta \sigma)|^2 \idx} \\
	& ~~~~ + 2 \sum_{\partial\in\lbrace \ddh, \dz \rbrace} \int (\partial v \cdot \nablah +  \partial w \dz )\bigl(\dt \partial \sigma - \Delta \partial \sigma\bigr) \Delta (\dt \sigma - \Delta \sigma)  \idx \\
	& ~~~~ + \int (\Delta v \cdot \nablah +  \Delta w \dz )\bigl(\dt \sigma - \Delta \sigma\bigr) \Delta (\dt \sigma - \Delta \sigma)  \idx. 
	\end{align*}
	Thus, taking the $ L^2 $-inner product of \eqref{eq:mixed} with $ \Delta^2 (\dt \sigma - \Delta \sigma) $ yields, after substituting the identities above, that
	\begin{equation}\label{est-204}
	\begin{aligned}
	& \dfrac{1}{2}\dfrac{d}{dt} \norm{\nabla^2 (\dt \sigma - \Delta \sigma), \nabla^2 \nablah \sigma, \dfrac{\nabla^2 \dz \sigma}{\varepsilon} }{\Lnorm{2}}^2 + \norm{\nabla^3 \nablah \sigma, \dfrac{\nabla^3 \dz \sigma}{\varepsilon}}{\Lnorm{2}}^2 \\
	& ~~~~ = \sum_{j=1}^4 \int \Delta \mathcal J_j \times \Delta (\dt\sigma - \Delta \sigma)\idx + \dfrac{1}{2} \int  (\dvh v + \dz w ) |\Delta (\dt \sigma - \Delta \sigma)|^2 \idx\\
	& ~~~~  - 2 \sum_{\partial\in\lbrace \ddh, \dz \rbrace} \int (\partial v \cdot \nablah +  \partial w \dz )\bigl(\dt \partial \sigma - \Delta \partial \sigma\bigr) \Delta (\dt \sigma - \Delta \sigma)  \idx \\
	& ~~~~ - \int (\Delta v \cdot \nablah +  \Delta w \dz )\bigl(\dt \sigma - \Delta \sigma\bigr) \Delta (\dt \sigma - \Delta \sigma)  \idx\\
	& ~~~~ \lesssim \norm{\nabla v,\nabla w}{\Lnorm{\infty}}  \norm{\nabla^2 (\dt \sigma - \Delta \sigma)}{\Lnorm{2}}^2 +  \norm{\nabla^2 v,\nabla^2 w}{\Lnorm{3}} \norm{\nabla (\dt \sigma - \Delta \sigma)}{\Lnorm{6}} \\
	& ~~~~ \times \norm{\nabla^2 (\dt \sigma - \Delta \sigma)}{\Lnorm{2}} + \sum_{j=1}^4 \norm{\nabla^2 \mathcal J_j}{\Lnorm{2}} \norm{\nabla^2(\dt \sigma - \Delta \sigma)}{\Lnorm{2}} \\
	& ~~~~ \lesssim \norm{v,w}{\Hnorm{3}} \norm{ \dt \sigma - \Delta \sigma}{\Hnorm{2}}^2 + \sum_{j=1}^4 \norm{\mathcal J_j}{\Hnorm{2}} \norm{\dt \sigma - \Delta \sigma}{\Hnorm{2}},
	\end{aligned}
	\end{equation}
	where $ \mathcal J_j $, $ j = 1,2,3,4$, are as in \eqref{id:nonlinearity-2}.

	On the other hand, we apply integration by parts to obtain the following integral identities:
	\begin{align*}
	& \int \dt( \dt \sigma - \Delta \sigma ) \Delta^2 \dt \sigma \idx  = \dfrac{1}{2} \dfrac{d}{dt} \norm{\nabla^2 \dt \sigma}{\Lnorm{2}}^2 + \norm{\nabla^3 \dt \sigma}{\Lnorm{2}}^2, \\
	& \int  ( -\deltah \sigma - \dfrac{1}{\varepsilon^2}\partial_{zz}\sigma ) \Delta^2 \dt \sigma \idx = \dfrac{1}{2} \dfrac{d}{dt} \norm{ \nabla^2 \nablah \sigma, \dfrac{\nabla
			^2 \dz \sigma}{\varepsilon} }{\Lnorm{2}}^2, \\
	& \int (v \cdot \nablah + w \dz )\bigl(\dt \sigma - \Delta \sigma\bigr) \Delta^2 \dt \sigma \idx \\
	& ~~~~  =  - \sum_{\partial\in\lbrace \ddh, \dz \rbrace}\int ( \partial v \cdot \nablah + \partial w\dz) (\dt \sigma -\Delta \sigma)  \partial \Delta \dt \sigma \idx \\
	& ~~~~ - \sum_{\partial\in\lbrace \ddh, \dz \rbrace}\int (  v \cdot \nablah +  w\dz)  (\dt \partial \sigma -\Delta \partial \sigma) \partial \Delta \dt \sigma \idx.
	\end{align*}
	Therefore, taking the $ L^2 $-inner product of \eqref{eq:mixed} with $ \Delta^2 \dt \sigma $ yields, after substituting the identities above, that 
	\begin{equation}\label{est-205}
	\begin{aligned}
	& \dfrac{1}{2} \dfrac{d}{dt} \norm{ \nabla^2 \dt\sigma, \nabla^2 \nablah \sigma, \dfrac{\nabla^2 \dz \sigma}{\varepsilon} }{\Lnorm{2}}^2 + \norm{\nabla^3 \dt \sigma}{\Lnorm{2}}^2\\
	& = \sum_{j = 1}^4 \int \Delta \mathcal J_j \times \Delta \dt\sigma \idx  
	 + \sum_{\partial\in\lbrace \ddh, \dz \rbrace} \int (  v \cdot \nablah +  w\dz) \partial (\dt \sigma -\Delta \sigma) \partial \Delta \dt \sigma \idx\\
	& ~~~~ + \sum_{\partial\in\lbrace \ddh, \dz \rbrace} \int ( \partial v \cdot \nablah + \partial w\dz) (\dt \sigma -\Delta \sigma) \partial \Delta \dt \sigma \idx \\
	& ~~~~ \lesssim \norm{\nabla v, \nabla  w}{\Lnorm{3}} \norm{\nabla(\dt\sigma - \Delta \sigma)}{\Lnorm{6}} \norm{\nabla \Delta \dt \sigma}{\Lnorm{2}} \\
	& ~~~~ + \norm{ v,   w}{\Lnorm{\infty}} \norm{\nabla^2(\dt\sigma - \Delta \sigma)}{\Lnorm{2}} \norm{\nabla \Delta \dt \sigma}{\Lnorm{2}} + \sum_{j=1}^4 \norm{\nabla^2 \mathcal J_j}{\Lnorm{2}} \norm{\nabla^2\dt \sigma }{\Lnorm{2}} \\
	& ~~~~ \lesssim \norm{v,w}{\Hnorm{2}} \norm{ \dt \sigma - \Delta \sigma}{\Hnorm{2}}\norm{\dt\sigma}{\Hnorm{3}} + \sum_{j=1}^4 \norm{\mathcal J_j}{\Hnorm{2}} \norm{\dt \sigma}{\Hnorm{2}},
	\end{aligned}
	\end{equation}
	where $ \mathcal J_j $, $ j = 1,2,3,4$, are as in \eqref{id:nonlinearity-2}.

	\subsection{Uniform-in-$\varepsilon$ estimates for the nonlinearities}\label{subsec:est-nonlinearities}
	In this section, we will estimate $ \norm{\mathcal I_i}{\Hnorm{2}}, \norm{\mathcal J_j}{\Hnorm{2}} $, $ i = 1,2, j = 1,2, 3, 4 $. 
	
	\begin{prop}\label{lm:nonlinearity}
	$ I_i, J_j $, $ i = 1,2, ~ j = 1,2,3,4$, given by \eqref{id:nonlinearity-1} and \eqref{id:nonlinearity-2}, satisfy the following estimates:
		\begin{gather}
		\label{est-nonlinearity-i-1}
		\norm{\mathcal I_i}{\Hnorm{2}} \leq C E^2, ~~~~ i = 1,2, \\
		\label{est-nonlinearity-j}
		\norm{\mathcal J_j}{\Hnorm{2}} \leq C E  ( E + D), ~~~~ j = 1,2,3,4,
		\end{gather}
		where $E$ and $D$ are as in \eqref{energy} and \eqref{dissipation}, respectively, and $ C \in (0,\infty) $ is some generic constant, independent of $ \varepsilon $.
	\end{prop}
	
	\begin{proof}
		We will only sketch the estimates for 
		\begin{equation*}
		\norm{\nabla^2 \mathcal I_i, \nabla^2 \mathcal J_j}{\Lnorm{2}}, ~~ i = 1,2, ~ j = 1,2,3,4.
		\end{equation*}
		The rests of \eqref{est-nonlinearity-i-1} and \eqref{est-nonlinearity-j}  can be derived via similar arguments. We organize the proof in the following order:
		\begin{equation*}
		\nabla^2\mathcal J_4, \nabla^2\mathcal J_3, \nabla^2\mathcal J_1, \nabla^2\mathcal J_2, \nabla^2\mathcal I_{1}, \nabla^2\mathcal I_{2}.
		\end{equation*}
		
		{\par\noindent\bf Estimate for $ \norm{\nabla^2 \mathcal J_4}{\Lnorm{2}} $: }
		From the definition, one can write
		\begin{equation*}
		\mathcal J_4 = \underbrace{\dz(v \cdot \nablah w)}_{\mathcal J_{41}} + \underbrace{\dz(w \dz w)}_{\mathcal J_{42}}.
		\end{equation*}
		Then applying the H\"older and Sobolev embedding inequalities yields
		\begin{align*}
		& \norm{\nabla^2 \mathcal J_{41}}{\Lnorm{2}} \lesssim \sum_{\partial \in \lbrace \ddh, \dz \rbrace} \biggl( \norm{\partial^2 \dz v \cdot \nablah w}{\Lnorm{2}} +  \norm{\partial \dz v \cdot \nablah \partial w}{\Lnorm{2}} \\
		& ~~~~ ~~~~ + \norm{ \dz v \cdot \nablah \partial^2 w}{\Lnorm{2}} + \norm{v \cdot \nablah \dz \partial^2 w}{\Lnorm{2}} + \norm{\partial v \cdot \nablah \dz \partial w}{\Lnorm{2}}\\
		& ~~~~ ~~~~ + \norm{\partial^2 v \cdot \nablah \dz w}{\Lnorm{2}}
		\biggr)  
		\lesssim  \sum_{\partial \in \lbrace \ddh, \dz \rbrace} \biggl( \norm{\partial^2 \dz v}{\Lnorm{2}}\norm{\nablah w}{\Lnorm{\infty}}\\
		& ~~~~ ~~~~ +  \norm{\partial \dz v}{\Lnorm{6}} \norm{\nablah \partial w}{\Lnorm{3}} + \norm{ \dz v}{\Lnorm{\infty}} \norm{\nablah \partial^2 w}{\Lnorm{2}} \\
		& ~~~~ ~~~~ + \norm{v}{\Lnorm{\infty}}\norm{\nablah \dz \partial^2 w}{\Lnorm{2}} + \norm{\partial v}{\Lnorm{3}}\norm{\nablah \dz \partial w}{\Lnorm{6}}\\
		& ~~~~ ~~~~ + \norm{\partial^2 v}{\Lnorm{2}}\norm{\nablah \dz w}{\Lnorm{\infty}}
		\biggr) \lesssim \norm{v}{\Hnorm{3}}\norm{w}{\Hnorm{3}} + \norm{v}{\Hnorm{2}}\norm{\dz w}{\Hnorm{3}}  \\
		& ~~~~ \lesssim E D,  \\
		& \norm{\nabla^2 \mathcal J_{42}}{\Lnorm{2}} \lesssim \sum_{\partial \in \lbrace \ddh, \dz \rbrace} \biggl( \norm{\partial^2 \dz w \dz w}{\Lnorm{2}} + \norm{\partial \dz w \dz \partial w }{\Lnorm{2}} \\
		& ~~~~ ~~~~ + \norm{\dz w \dz \partial^2 w}{\Lnorm{2}}+ \norm{\partial^2 w \partial_{zz} w }{\Lnorm{2}} + \norm{\partial w \partial_{zz} \partial w }{\Lnorm{2}}\\
		& ~~~~ ~~~~ + \norm{w\partial_{zz} \partial^2 w}{\Lnorm{2}}  \biggr) 
		\lesssim \sum_{\partial \in \lbrace \ddh, \dz \rbrace} \biggl( \norm{\partial^2 \dz w}{\Lnorm{2}}\norm{\dz w}{\Lnorm{\infty}} \\
		& ~~~~ ~~~~ + \norm{\partial \dz w}{\Lnorm{3}}\norm{\dz \partial w }{\Lnorm{6}} 
		+ \norm{\dz w}{\Lnorm{\infty}} \norm{\dz \partial^2 w}{\Lnorm{2}}\\
		& ~~~~ ~~~~ + \norm{\partial^2 w}{\Lnorm{2}} \norm{\partial_{zz} w }{\Lnorm{\infty}} + \norm{\partial w}{\Lnorm{3}} \norm{\partial_{zz} \partial w }{\Lnorm{6}}\\
		& ~~~~ ~~~~ + \norm{w}{\Lnorm{\infty}}\norm{\partial_{zz} \partial^2 w}{\Lnorm{2}}  \biggr) \lesssim \norm{\dz w}{\Hnorm{2}} \norm{\dz w}{\Hnorm{2}} + \norm{w}{\Hnorm{2}} \norm{\dz w}{\Hnorm{3}} \\
		& ~~~~ \lesssim E^2 + E  D.
		\end{align*}
		
		{\par\noindent\bf Estimate for $ \norm{\nabla^2 \mathcal J_3}{\Lnorm{2}} $: } We write
		\begin{equation*}
		\mathcal J_3 = \underbrace{\dvh ( v\cdot \nablah v )}_{\mathcal J_{31}} + \underbrace{\dvh (w \dz v)}_{\mathcal J_{32}}.
		\end{equation*}
		Then applying the H\"older and Sobolev embedding inequalities yields
		\begin{align*}
		& \norm{\nabla^2 \mathcal J_{31}}{\Lnorm{2}} \lesssim \sum_{\partial \in \lbrace \ddh, \dz \rbrace} \biggl(\norm{\partial^3 v \cdot \nablah v}{\Lnorm{2}} + \norm{\partial^2 v \cdot \nablah \partial v}{\Lnorm{2}} \\
		& ~~~~ ~~~~ + \norm{\partial v \cdot \nablah \partial^2 v}{\Lnorm{2}} + \norm{v \cdot \nablah \partial^3 v}{\Lnorm{2}} 
		\biggr) \\
		& ~~~~ \lesssim \sum_{\partial \in \lbrace \ddh, \dz \rbrace} \biggl(\norm{\partial^3 v}{\Lnorm{2}} \norm{\nablah v}{\Lnorm{\infty}} + \norm{\partial^2 v}{\Lnorm{6}} \norm{\nablah \partial v}{\Lnorm{3}} \\
		& ~~~~ ~~~~ + \norm{\partial v}{\Lnorm{\infty}} \norm{\nablah \partial^2 v}{\Lnorm{2}} + \norm{v}{\Lnorm{\infty}} \norm{\nablah \partial^3 v}{\Lnorm{2}} 
		\biggr) \\
		& ~~~~ \lesssim \norm{v}{\Hnorm{3}}\norm{v}{\Hnorm{3}} + \norm{v}{\Hnorm{2}} \norm{v}{\Hnorm{4}} \lesssim E^2 + E D, \\
		& \norm{\nabla^2 \mathcal J_{32}}{\Lnorm{2}} \lesssim \norm{w}{\Hnorm{3}}\norm{v}{\Hnorm{3}} + \norm{w}{\Hnorm{2}}\norm{v}{\Hnorm{4}}  
		 \lesssim E D.
		\end{align*}

		{\par\noindent\bf Estimate for $ \norm{\nabla^2 \mathcal J_1}{\Lnorm{2}} $: }
		\begin{equation*}
		\mathcal J_1 = \underbrace{\dt v \cdot \nablah \sigma}_{\mathcal J_{11}} - \underbrace{\Delta v \cdot \nablah \sigma}_{\mathcal J_{12}} -2  \underbrace{\nabla v : \nablah \nabla \sigma}_{\mathcal J_{13}}.
		\end{equation*}
		We list the estimates in the following:
		\begin{align*}
		& \norm{\nabla^2 \mathcal J_{11}}{\Lnorm{2}} \lesssim \sum_{\partial \in \lbrace \ddh, \dz \rbrace} \biggl( \norm{\partial^2 \dt v\cdot\nablah \sigma}{\Lnorm{2}} + \norm{\partial \dt v\cdot \nablah \partial \sigma}{\Lnorm{2}} \\
		& ~~~~ ~~~~ + \norm{\dt v \cdot \nablah \partial^2 \sigma}{\Lnorm{2}} \biggr) \lesssim \sum_{\partial \in \lbrace \ddh, \dz \rbrace} \biggl( \norm{\partial^2 \dt v}{\Lnorm{2}} \norm{\nablah \sigma}{\Lnorm{\infty}} \\
		& ~~~~ ~~~~ + \norm{\partial \dt v}{\Lnorm{3}} \norm{\nablah \partial \sigma}{\Lnorm{6}} + \norm{\dt v}{\Lnorm{\infty}} \norm{\nablah \partial^2 \sigma}{\Lnorm{2}} 
		\biggr) \\
		& ~~~~ \lesssim \norm{\dt v}{\Hnorm{2}} \norm{\sigma}{\Hnorm{3}} \lesssim E D, \\
		& \norm{\nabla^2 \mathcal J_{12}}{\Lnorm{2}} \lesssim \norm{\Delta v}{\Hnorm{2}} \norm{\sigma}{\Hnorm{3}} \lesssim E D, \\
		& \norm{\nabla^2 \mathcal J_{13}}{\Lnorm{2}} \lesssim \sum_{\partial \in \lbrace \ddh, \dz \rbrace }  \biggl(\norm{\partial^2 \nabla v : \nablah \nabla\sigma}{\Lnorm{2}}  + \norm{\partial\nabla v : \nablah \nabla \partial \sigma}{\Lnorm{2}} \\
		& ~~~~ ~~~~ + \norm{\nabla v : \nablah \nabla\partial^2 \sigma}{\Lnorm{2}}  
		\biggr) \lesssim \sum_{\partial \in \lbrace \ddh, \dz \rbrace }  \biggl(\norm{\partial^2 \nabla v}{\Lnorm{2}} \norm{\nablah \nabla\sigma}{\Lnorm{\infty}} \\
		& ~~~~ ~~~~ + \norm{\partial\nabla v}{\Lnorm{3}} \norm{\nablah \nabla \partial \sigma}{\Lnorm{6}} 
		+ \norm{\nabla v}{\Lnorm{\infty}} \norm{\nablah \nabla\partial^2 \sigma}{\Lnorm{2}}  
		\biggr) \\
		& ~~~~ \lesssim \norm{v}{\Hnorm{3}} \norm{\sigma}{\Hnorm{4}} \lesssim E^2.
		\end{align*}
		
		{\par\noindent\bf Estimate for $ \norm{\nabla^2 \mathcal J_2}{\Lnorm{2}} $: } We write $ \mathcal J_2 $ as
		\begin{equation*}
		\mathcal J_2 = \underbrace{\dt (\varepsilon w) \dfrac{\dz \sigma}{\varepsilon}}_{\mathcal J_{21}}  - \underbrace{\Delta(\varepsilon w) \dfrac{\dz \sigma}{\varepsilon}}_{\mathcal J_{22}} - 2 \underbrace{\nabla (\varepsilon w) \cdot \dfrac{\nabla \dz \sigma}{\varepsilon}}_{\mathcal J_{23}}. 
		\end{equation*}
		Then applying the H\"older and Sobolev embedding inequalities implies,
		\begin{align*}
		& \norm{\nabla^2 \mathcal J_{21}}{\Lnorm{2}} \lesssim \sum_{\partial \in \lbrace \ddh, \dz \rbrace} \biggl( \norm{\partial^2 \dt(\varepsilon w) \dfrac{\dz \sigma}{\varepsilon}}{\Lnorm{2}} + \norm{\partial \dt(\varepsilon w) \dfrac{\dz \partial \sigma}{\varepsilon}}{\Lnorm{2}} \\
		& ~~~~ ~~~~ + \norm{\dt(\varepsilon w) \dfrac{\dz \partial^2 \sigma}{\varepsilon}}{\Lnorm{2}} \biggr) \lesssim \sum_{\partial \in \lbrace \ddh, \dz \rbrace} \biggl( \norm{\partial^2 \dt(\varepsilon w)}{\Lnorm{2}} \norm{\dfrac{\dz \sigma}{\varepsilon}}{\Lnorm{\infty}} \\
		& ~~~~ ~~~~ + \norm{\partial \dt(\varepsilon w)}{\Lnorm{6}} \norm{\dfrac{\dz \partial \sigma}{\varepsilon}}{\Lnorm{3}} 
		+ \norm{\dt(\varepsilon w)}{\Lnorm{\infty}} \norm{\dfrac{\dz \partial^2 \sigma}{\varepsilon}}{\Lnorm{2}} \biggr)  \\
		& ~~~~ \lesssim \norm{\dt (\varepsilon w)}{\Hnorm{2}} \norm{\dfrac{\dz \sigma}{\varepsilon}}{\Hnorm{2}} \lesssim E D, \\
		& \norm{\nabla^2 \mathcal J_{22}}{\Lnorm{2}} \lesssim \norm{\Delta (\varepsilon w)}{\Hnorm{2}} \norm{\dfrac{\dz \sigma}{\varepsilon}}{\Hnorm{2}} \lesssim  E D, \\
		& \norm{\nabla^2 \mathcal J_{23}}{\Lnorm{2}} \lesssim \sum_{\partial\in \lbrace \ddh,\dz\rbrace} \biggl( \norm{\partial^2 \nabla(\varepsilon w) \cdot \dfrac{\nabla \dz \sigma}{\varepsilon}}{\Lnorm{2}} + \norm{\partial \nabla(\varepsilon w) \cdot \dfrac{\nabla \dz \partial \sigma}{\varepsilon}}{\Lnorm{2}} \\
		& ~~~~ ~~~~ + \norm{\nabla(\varepsilon w) \cdot \dfrac{\nabla \dz \partial^2 \sigma}{\varepsilon}}{\Lnorm{2}} \lesssim 
		\biggr) 
		\lesssim \sum_{\partial\in \lbrace \ddh,\dz\rbrace} \biggl( \norm{\partial^2 \nabla(\varepsilon w)}{\Lnorm{2}} \norm{ \dfrac{\nabla \dz \sigma}{\varepsilon}}{\Lnorm{\infty}} \\
		& ~~~~ ~~~~ + \norm{\partial \nabla(\varepsilon w)}{\Lnorm{6}} \norm{\dfrac{\nabla \dz \partial \sigma}{\varepsilon}}{\Lnorm{3}} 
		+ \norm{\nabla(\varepsilon w)}{\Lnorm{\infty}} \norm{\dfrac{\nabla \dz \partial^2 \sigma}{\varepsilon}}{\Lnorm{2}} \biggr) \\
		& ~~~~ \lesssim \norm{\varepsilon w}{\Hnorm{3}} \norm{\dfrac{\dz \sigma}{\varepsilon}}{\Hnorm{3}} \lesssim E D.
		\end{align*}
		
		{\par\noindent\bf Estimate for $ \norm{\nabla^2 \mathcal I_1}{\Lnorm{2}}, \norm{\nabla^2 \mathcal I_2}{\Lnorm{2}} $: } 
		The estimates are similar. They are,
		\begin{align*}
		& \norm{ \nabla^2 \mathcal I_1}{\Lnorm{2}} \lesssim \sum_{\partial \in \lbrace \ddh, \dz \rbrace} \biggl( \norm{\partial^2 v \cdot \nablah (v, \varepsilon w)}{\Lnorm{2}} + \norm{\partial v \cdot \nablah \partial( v, \varepsilon w)}{\Lnorm{2}} \\
		& ~~~~ ~~~~ + \norm{v \cdot \nablah \partial^2 (v, \varepsilon w)}{\Lnorm{2}} \biggr) \lesssim \sum_{\partial \in \lbrace \ddh, \dz \rbrace} \biggl( \norm{\partial^2 v}{\Lnorm{2}} \norm{\nablah (v, \varepsilon w)}{\Lnorm{\infty}} \\
		& ~~~~ ~~~~ + \norm{\partial v}{\Lnorm{3}} \norm{\nablah \partial( v, \varepsilon w)}{\Lnorm{6}} 
		+ \norm{v}{\Lnorm{\infty}} \norm{\nablah \partial^2 (v, \varepsilon w)}{\Lnorm{2}} \biggr) \\
		& ~~~~ \lesssim \norm{v}{\Hnorm{2}} \norm{v, \varepsilon w}{\Hnorm{3}} \lesssim E^2, \\
		& \norm{\nabla^2 \mathcal I_2}{\Lnorm{2}} \lesssim \sum_{\partial \in \lbrace \ddh, \dz \rbrace} \biggl( \norm{\partial^2 w \dz (v, \varepsilon w)}{\Lnorm{2}} + \norm{\partial w \dz \partial( v, \varepsilon w)}{\Lnorm{2}} \\
		& ~~~~ ~~~~ + \norm{w \dz \partial^2 (v, \varepsilon w)}{\Lnorm{2}} \biggr) \lesssim \sum_{\partial \in \lbrace \ddh, \dz \rbrace} \biggl( \norm{\partial^2 w}{\Lnorm{2}} \norm{\dz(v, \varepsilon w)}{\Lnorm{\infty}} \\
		& ~~~~ ~~~~ + \norm{\partial w}{\Lnorm{3}} \norm{\dz \partial( v, \varepsilon w)}{\Lnorm{6}} 
		+ \norm{w}{\Lnorm{\infty}} \norm{\dz \partial^2 (v, \varepsilon w)}{\Lnorm{2}} \biggr) \\
		& ~~~~ \lesssim \norm{w}{\Hnorm{2}} \norm{v, \varepsilon w}{\Hnorm{3}} \lesssim E^2.
		\end{align*}
	\end{proof}

	\subsection{Summary of uniform-in-$\varepsilon$ estimates}\label{subsec:uniform-est}
	To summarize the estimates in the previous sections, let
	\begin{align}
	& \label{instant-energy}
	\begin{aligned}
	&	\mathcal E :=  \norm{v, \nabla v, \varepsilon w, \nabla (\varepsilon w), \sigma, \dt \sigma - \Delta \sigma, \dt \sigma}{\Lnorm{2}}^2 +2 \norm{\nablah \sigma, \dfrac{\dz \sigma}{\varepsilon}}{\Lnorm{2}}^2 \\
	& ~~ + \norm{ \nabla^3 v, \nabla^3(\varepsilon w),
		\nabla^2 \dt \sigma ,\nabla^2(\dt \sigma - \Delta \sigma) }{\Lnorm{2}}^2 + 2 \norm{ \nabla^2 \nablah \sigma, \dfrac{\nabla^2 \dz \sigma}{\varepsilon}}{\Lnorm{2}}^2,
	\end{aligned} \\
	& \label{instant-dissipation}
	\begin{aligned}
	& \mathcal D := \norm{\dt v, \nabla^2 v, \dt(\varepsilon w), \nabla^2 ( \varepsilon w)}{\Lnorm{2}}^2 + 2 \norm{\nabla \dt \sigma ,\nabla \nablah \sigma, \dfrac{\nabla\dz \sigma}{\varepsilon}}{\Lnorm{2}}^2 \\
	& ~~ + \norm{ \nabla^2 \dt v, \nabla^4 v, \nabla^2 \dt (\varepsilon w), \nabla^4 (\varepsilon w)}{\Lnorm{2}}^2 + 2 \norm{
		\nabla^3 \dt \sigma , \nabla^3 \nablah \sigma, \dfrac{\nabla^3 \dz \sigma}{\varepsilon}}{\Lnorm{2}}^2.
	\end{aligned}
	\end{align}
	Then, from the estimates in \eqref{est-w-1}, \eqref{est-w-2}, the definitions in \eqref{energy}, \eqref{dissipation}, \eqref{instant-energy}, and \eqref{instant-dissipation}, one concludes, after applying the triangle inequality and the interpolation inequality for Sobolev spaces, that
	\begin{equation}\label{eqvl-functionals}
	\begin{gathered}
	\mathcal E \leq E^2 \leq \mathfrak N(\mathcal E), \\
	D \leq \mathfrak N(\mathcal E,1) \sqrt{\mathcal D} + \mathfrak N (\mathcal E). 
	\end{gathered}
	\end{equation}
	Let us recall once again that $ \mathfrak N(\cdot) $ is a locally Lipschitz in its argument. 
	
	Applying the Cauchy--Schwarz inequality implies that
	\begin{equation}\label{est-301}
	\dfrac{d}{dt} \norm{v, \varepsilon w,\sigma}{\Lnorm{2}}^2 \leq 2 \norm{v, \varepsilon w,\sigma}{\Lnorm{2}}\norm{\dt v, \dt(\varepsilon w),\dt \sigma}{\Lnorm{2}}. 
	\end{equation}
	Thus, \eqref{est-003}, \eqref{est-density-001}, \eqref{est-203}, \eqref{est-204}, \eqref{est-205}, and \eqref{est-301} yield
	\begin{equation}\label{ineq:total-0}
	\begin{aligned}
	& \dfrac{d}{dt} \mathcal E + \mathcal D \leq \mathfrak N ( E )  +  C (E + E^2) D + C \sum_{i = 1}^2 \norm{\mathcal I_{i}}{\Hnorm{2}}^2 + C E \sum_{j=1}^4  \norm{\mathcal J_j}{\Hnorm{2}} \\
	& ~~~~ \leq \mathfrak N(E) + \mathfrak N(E) D \leq \mathfrak N (\mathcal E) + \mathfrak N(\mathcal E,1) \sqrt{\mathcal D},
	\end{aligned}
	\end{equation}
	where the second and the last inequalities are consequences of \eqref{est-nonlinearity-i-1}, \eqref{est-nonlinearity-j}, and \eqref{eqvl-functionals}. 
	Therefore, \eqref{ineq:total-0} implies, after applying the Cauchy--Schwarz inequality, that
	\begin{equation}\label{ineq:total-1}
	\dfrac{d}{dt} \mathcal E + \dfrac{1}{2} \mathcal D \leq \mathfrak N(\mathcal E,1). 
	\end{equation}
	Then for some $ T^* \in (0,\infty) $, independent of $ \varepsilon $, provided  that $ \mathcal E(0) $ is finite, we have
	\begin{equation}\label{uniform-est-1}
	\sup_{0\leq t \leq T^* } \mathcal E(t) + \int_0^{T^*} \mathcal D(t) \,dt < \infty. 
	\end{equation}
	Together with \eqref{eqvl-functionals}, we have
	\begin{equation}\label{uniform-est-2}
	\sup_{0\leq t \leq T^* } E^2(t) + \int_0^{T^*} D^2(t)\,dt \leq C < \infty, 
	\end{equation}
	where $ C \in (0,\infty) $ depends only on the initial data and is independent of $ \varepsilon $.

\section{Convergence to the compressible primitive equations}\label{sec:convergence-cpe}

The definitions of $ E, D $ in \eqref{energy}, \eqref{dissipation}, respectively, and estimate
\eqref{uniform-est-2} imply that there are a time $ T^* \in (0,\infty) $ and a constant $ C \in (0,\infty) $, which are independent of $ \varepsilon $, such that 
\begin{equation}\label{uniform-est-con-ill-pre}
\begin{gathered}
 \norm{v_\varepsilon, \varepsilon w_\varepsilon}{L^\infty(0, T^{*};\Hnorm{3} )} + \norm{\dt v_\varepsilon, \dt (\varepsilon w_\varepsilon)}{L^2(0, T^{*};\Hnorm{2} )} +
\norm{v_\varepsilon, \varepsilon w_\varepsilon}{L^2(0, T^{*};\Hnorm{4} )} \\
 + \norm{\dt \sigma_\varepsilon, \nablah \sigma_\varepsilon, \dfrac{\dz \sigma_\varepsilon}{\varepsilon}}{L^\infty(0,T^{*};H^2)} 
+ \norm{\sigma_\varepsilon}{L^\infty(0,T^{*};\Hnorm{4})} \\
 + \norm{\dt \sigma_\varepsilon, \nablah \sigma_\varepsilon, \dfrac{\dz \sigma_\varepsilon}{\varepsilon}}{L^2(0,T^{*};H^3)} 
 + \norm{w_\varepsilon,\dz w_\varepsilon}{L^\infty(0,T^{*};H^2)} \\
 + \norm{w_\varepsilon,\dz w_\varepsilon}{L^2(0,T^{*};H^3)} < C.
\end{gathered}
\end{equation}
Then from \eqref{uniform-est-con-ill-pre}, one can conclude that there exist
\begin{equation}\label{con-ill-pre-functions}
\begin{gathered}
\sigma^* \in L^\infty(0,T^{*};H^4), ~~ \dt \sigma^* \in L^\infty(0,T^*;H^2) \cap L^2(0,T^*;H^3), \\
v^*\in L^\infty(0,T^*;H^3) \cap L^2(0,T^*;H^4), ~~ \dt v^* \in L^2(0,T^*;H^2), \\
w^*, \dz w^* \in L^\infty(0,T^*;H^2)\cap L^2(0,T^*;H^3),
\end{gathered}
\end{equation}
such that for a subsequence of $ \lbrace (\sigma_\varepsilon, v_\varepsilon, w_\varepsilon) \rbrace $, as $ \varepsilon \rightarrow 0^+ $, one has, 
\begin{align*}
 \sigma_\varepsilon & \buildrel\ast\over\rightharpoonup \sigma^* & ~ \text{weak-$\ast$ in} ~ &  L^\infty(0,T^*;H^4), \\
 \sigma_\varepsilon & \rightarrow \sigma^* & ~ \text{in} ~  & L^\infty(0,T^*;H^3)\cap C([0,T^*];H^3),\\
 \dt \sigma_\varepsilon, w_\varepsilon, \dz w_\varepsilon & \buildrel\ast\over\rightharpoonup \dt \sigma^*, w^*, \dz w^* & ~ \text{weak-$\ast$ in} ~ & L^\infty(0,T^*;H^2), \\
 \dt \sigma_\varepsilon , w_\varepsilon, \dz w_\varepsilon & \rightharpoonup \dt \sigma^*, w^*, \dz w^* & ~ \text{weakly in} ~ & L^2(0,T^*;H^3), \\
 v_\varepsilon & \buildrel\ast\over\rightharpoonup v^* & ~ \text{weak-$\ast$ in} ~ & L^\infty(0,T^*;H^3), \\
  v_\varepsilon & \rightarrow v^* & ~ \text{in} ~ & L^\infty(0,T^*;H^2)\cap C([0,T^*];H^2), \\
   v_\varepsilon & \rightharpoonup v^* & ~ \text{weakly in} ~ & L^2(0,T^*;H^4), \\
     \dt v_\varepsilon & \rightharpoonup \dt v^* & ~ \text{weakly in} ~ & L^2(0,T^*;H^2),
\end{align*}
where we have applied the Aubin compactness Lemma (see, e.g., \cite{Temam1984}). 
Then it is easy to verify that $ (\sigma^*, v^*, w^*) $ satisfies the compressible primitive equations \eqref{CPE}, after applying the above weak and strong convergences in each terms in system \eqref{CF}. We omit the details.

\section{The rates of convergence in terms of $ \varepsilon$}\label{sec:converging-rate}

In this section, we aim at investigating the rates of convergence of $ \dlsigma, \dlv, \dlv $ to $ 0 $, as $ \varepsilon \rightarrow 0^+ $. In section \ref{subsec:uniform-est-2}, we derive some additional uniform-in-$\varepsilon$ estimates, which will be used in section \ref{subsec:converging-rate} to derive estimates for the rates of convergence. In section \ref{subsec:proof-est-t}, we provide the proofs of some Propositions, which have been used in section \ref{subsec:uniform-est-2}. 

We shorten the notations by dropping the subscript $ \varepsilon $ in $ (\sigma_\varepsilon, v_\varepsilon, w_\varepsilon) $, i.e., $ (\sigma, v, w) = (\sigma_\varepsilon, v_\varepsilon, w_\varepsilon) $.

\subsection{Additional uniform-in-$\varepsilon$ bounds}\label{subsec:uniform-est-2}

The additional uniform-in-$\varepsilon$ estimates, in order to obtain the convergence rates, are basically the temporal derivative version of the estimates in section \ref{sec:uniform-est}. 
Denote by
\begin{align}
\label{instant-energy-t} & 
\begin{aligned}
& \mathcal E_1 := \norm{\dt v, \nabla \dt v, \dt (\varepsilon w), \nabla\dt (\varepsilon w ),\partial_t^2 \sigma - \Delta \dt \sigma, \dt^2 \sigma }{\Lnorm{2}}^2 \\
& ~~~~ + 2 \norm{\nablah \dt \sigma, \dfrac{\dz\dt\sigma}{\varepsilon}}{\Lnorm{2}}^2,
\end{aligned} \\
\label{instant-dissipation-t} &
\begin{aligned}
& \mathcal D_1 := \norm{\partial_t^2 v, \nabla^2 \dt v, \partial_t^2 (\varepsilon w), \nabla^2 \dt (\varepsilon w )}{\Lnorm{2}}^2 + 2 \norm{\nabla\nablah \dt \sigma, \dfrac{\nabla \dz \dt \sigma}{\varepsilon}}{\Lnorm{2}}^2 \\
& ~~~~ + 2 \norm{\nabla\partial_t^2 \sigma}{\Lnorm{2}}^2.
\end{aligned} 
\end{align}

After applying one temporal derivative to \subeqref{CF}{2}, \subeqref{CF}{3}, and \eqref{eq:mixed}, we obtain the following set of equations:
\begin{align}
\label{eq:v-t} & \begin{aligned}
& \partial_t v_t - \deltah v_t - \partial_{zz} v_t = - \dt \nablah \sigma - \dt ( v\cdot \nablah v ) - \dt (w \dz v ),
\end{aligned}
\\
\label{eq:w-t} & \begin{aligned}
& \dt w_t - \deltah w_t - \partial_{zz} w_t = - \dfrac{1}{\varepsilon^2} \dz \dt \sigma - \dt (v \cdot \nablah w ) - \dt ( w \dz w ),
\end{aligned}
\\
\label{eq:mixed-t} & \begin{aligned}
& \dt \bigl(\dt \sigma_t - \Delta \sigma_t \bigr) + (v \cdot\nablah + w \dz) \bigl( \dt \sigma_t - \Delta \sigma_t \bigr) - \deltah \sigma_t - \dfrac{1}{\varepsilon^2} \partial_{zz} \sigma_t \\
& ~~~~ 
= - (v_t \cdot\nablah + w_t \dz) \bigl( \dt \sigma - \Delta \sigma \bigr) + \sum_{j=1}^{4} (\mathcal{J}_j)_t,
\end{aligned}
\end{align}
where $ \mathcal J_j $, $ j = 1,2,3,4$, are as in \eqref{id:nonlinearity-2}.
Take the $ L^2 $-inner product of \eqref{eq:v-t}, \eqref{eq:w-t} with $ \dt v_t - \Delta v_t $, $ \varepsilon^2 ( \dt  w_t - \Delta w_t ) $, respectively. Then after applying integration by parts and the H\"older inequality in the resultant equations, similar arguments as in from \eqref{est-001} to \eqref{est-003} lead to
\begin{equation}\label{est-101}
\begin{aligned}
& \dfrac{d}{dt} \norm{\nabla \dt v, \varepsilon \nabla \dt w }{\Lnorm{2}}^2 + \norm{\partial_t^2 v, \nabla^2 \dt v, \varepsilon \partial_t^2 w, \varepsilon \nabla^2 \dt w}{\Lnorm{2}}^2 \\
& ~~~~ \leq C \norm{\nablah \dt \sigma, \dfrac{\dz\dt \sigma}{\varepsilon}}{\Lnorm{2}}^2 + C \sum_{i=1}^2 \norm{\dt \mathcal{I}_i}{\Lnorm{2}}^2,
\end{aligned}
\end{equation}
for some constant $ C \in (0,\infty) $, independent of $ \varepsilon $.
In addition, applying the Cauchy-Schwarz 
inequality yields,
\begin{equation}\label{est-1011}
\dfrac{d}{dt}\norm{\dt v, \dt (\varepsilon w)}{\Lnorm{2}}^2 \leq 2\norm{\dt v, \dt (\varepsilon w) }{\Lnorm{2}} \norm{\dt^2 v, \dt^2 (\varepsilon w )}{\Lnorm{2}}.
\end{equation}

On the other hand, taking the $ L^2 $-inner product of \eqref{eq:mixed-t} with $ \dt \sigma_t - \Delta \sigma_t $ yields, after applying integration by parts, the H\"older and Sobolev inequalities,
\begin{equation}\label{est-102}
\begin{aligned}
& \dfrac{1}{2} \dfrac{d}{dt} \norm{\partial_t^2 \sigma - \Delta \dt \sigma, \nablah \dt \sigma, \dfrac{\dz \dt \sigma}{\varepsilon}}{\Lnorm{2}}^2 + \norm{\nabla \nablah \dt\sigma, \dfrac{\nabla \dz \dt \sigma}{\varepsilon}}{\Lnorm{2}}^2 \\
& = \sum_{j=1}^4 \int \dt \mathcal J_j \times (\partial_t^2 \sigma - \Delta \dt \sigma )\idx + \dfrac{1}{2} \int (\dvh v + \dz w) |\partial_t^2 \sigma - \Delta \dt \sigma|^{2} \idx \\
& ~~~~ - \int (v_t \cdot\nablah + w_t \dz) \bigl( \dt \sigma - \Delta \sigma \bigr) \bigl(\partial_t^2 \sigma - \Delta \dt \sigma \bigr)\idx \\
& \lesssim \sum_{j=1}^4 \norm{\dt \mathcal J_j}{\Lnorm{2}}\norm{\partial_t^2 \sigma - \Delta \dt \sigma}{\Lnorm{2}} + \norm{\nablah v, \dz w}{\Lnorm{\infty}} \norm{\partial_t^2 \sigma - \Delta \dt \sigma}{\Lnorm{2}}^2 \\
& ~~~~ + {  \norm{v_t, w_t}{\Lnorm{3}} \norm{\nabla \dt \sigma, \nabla \Delta \sigma}{\Lnorm{6}} \norm{\partial_t^2 \sigma - \Delta \dt \sigma}{\Lnorm{2}}} \\
& \lesssim \sum_{j=1}^4 \norm{\dt \mathcal J_j}{\Lnorm{2}}\norm{\partial_t^2 \sigma - \Delta \dt \sigma}{\Lnorm{2}} + (\norm{v}{\Hnorm{3}} + \norm{\dz w}{\Hnorm{2}}) \norm{\partial_t^2 \sigma - \Delta \dt \sigma}{\Lnorm{2}}^2 \\
& ~~~~ + {  \norm{v_t, w_t}{\Hnorm{1}} \norm{\nabla \dt \sigma, \nabla \Delta \sigma}{\Hnorm{1}} \norm{\partial_t^2 \sigma - \Delta \dt \sigma}{\Lnorm{2}}}
. 
\end{aligned}
\end{equation}
In the meantime, after taking the $ L^2 $-inner product of \eqref{eq:mixed-t} with $ \partial_{t}^2 \sigma $ and applying integration by parts, the H\"older and Sobolev inequalities in the resultant equation, one can show that
\begin{equation}\label{est-103}
\begin{aligned}
& \dfrac{1}{2} \dfrac{d}{dt} \norm{\partial_t^2 \sigma, \nablah \dt \sigma, \dfrac{\dz \dt \sigma}{\varepsilon}}{\Lnorm{2}}^2 + \norm{\nabla \partial_t^2 \sigma }{\Lnorm{2}}^2 \\
& ~~~~ = \sum_{j=1}^4 \int \dt \mathcal J_j \times \partial_t^2 \sigma \idx
+ \dfrac{1}{2} \int (\dvh v + \dz w) |\partial_t^2 \sigma|^2 \idx \\
& - \int \biggl\lbrack (\dvh v + \dz w) \Delta \dt \sigma \biggr\rbrack  \partial_t^2 \sigma \idx - \int \biggl\lbrack (v \cdot \nablah + w \dz) \partial_t^2 \sigma \biggr\rbrack \Delta \dt \sigma \idx\\
& ~~~~ - \int \biggl\lbrack (v_t \cdot \nablah + w_t \dz )\bigl( \dt\sigma - \Delta \sigma \bigr) \biggr\rbrack \partial_t^2 \sigma \idx \\
& ~~~~ \lesssim \sum_{j=1}^4 \norm{\dt \mathcal J_j}{\Lnorm{2}} \norm{\partial_t^2 \sigma}{\Lnorm{2}} + \norm{\dvh v, \dz w}{\Lnorm{\infty}} \norm{\Delta \partial_t \sigma, \partial_t^2 \sigma}{\Lnorm{2}}^2\\
& ~~~~ + \norm{v, w}{\Lnorm{\infty}} \norm{\nabla\partial_t^2 \sigma}{\Lnorm{2}} \norm{\Delta\dt \sigma}{\Lnorm{2}} 
\\& ~~~~
+ { \norm{\dt v, \dt w}{\Lnorm{3}} \norm{\nabla \dt \sigma, \nabla \Delta\sigma }{\Lnorm{6}} \norm{\partial_t^2 \sigma}{\Lnorm{2}} 
} \\
& ~~~~ \lesssim \sum_{j=1}^4 \norm{\dt \mathcal J_j}{\Lnorm{2}} \norm{\partial_t^2 \sigma}{\Lnorm{2}} + (\norm{v}{\Hnorm{3}} + \norm{\dz w}{\Hnorm{2}}) \norm{\Delta \partial_t \sigma, \partial_t^2 \sigma}{\Lnorm{2}}^2\\
& ~~~~ + \norm{v, w}{\Hnorm{2}} \norm{\nabla\partial_t^2 \sigma}{\Lnorm{2}} \norm{\Delta\dt \sigma}{\Lnorm{2}} 
\\& ~~~~
+ { \norm{\dt v, \dt w}{\Hnorm{1}} \norm{\nabla \dt \sigma, \nabla \Delta\sigma }{\Hnorm{1}} \norm{\partial_t^2 \sigma}{\Lnorm{2}} 
} .
\end{aligned}
\end{equation}

Next, we state the estimates for $ \norm{\dt \mathcal I_i}{\Lnorm{2}}, \norm{\dt \mathcal J_j}{\Lnorm{2}} $, $, i = 1,2, j = 1,2,3,4 $ in the following: 
\begin{prop}\label{lm:nonlinearity-t} $ \mathcal I_i, \mathcal J_j $, $ i = 1,2, ~ j = 1,2,3,4$, given by \eqref{id:nonlinearity-1} and \eqref{id:nonlinearity-2}, satisfy the following estimates:
	\begin{gather}
	\norm{\dt \mathcal I_i}{\Lnorm{2}} \leq C  (E + E_1)^2, ~~~~ i = 1,2, \label{est-nonlinearity-i-t-1} \\
	\norm{\dt \mathcal J_j}{\Lnorm{2}} \leq C (E + E_1)  (E + E_1 +  D + D_1), ~~~~ j = 1,2,3,4, \label{est-nonlinearity-j-t}
	\end{gather}
	where $ C \in (0,\infty) $ is some generic constant, independent of $ \varepsilon $. Here, $ E, D, E_1, D_1 $ are defined in \eqref{energy}, \eqref{dissipation}, \eqref{energy-1}, and \eqref{dissipation-1}, respectively. 
\end{prop}

In addition, we state the estimates for $ \norm{\dt w, \dz \dt w}{\Lnorm{2}}, \norm{\dt w, \dz \dt w}{\Hnorm{1}} $ in the following: 
\begin{prop}\label{lm:vertical-velocity-t}
Let $ w $ be given as in \eqref{id:vertical-velocity-0}. Then it holds
	\begin{equation}\label{est-w-t}
	\begin{aligned}
	& \norm{\dt w, \dz \dt w}{\Lnorm{2}}\leq C e^{2\norm{\sigma}{\Hnorm{2}}} \bigl(\norm{\sigma}{\Hnorm{3}}^2 + \norm{\dt \sigma}{\Hnorm{1}}^2 + 1 \bigr) \bigl( \norm{\dt^2 \sigma}{\Lnorm{2}} \\
	& ~~~~ + \norm{\dt v}{\Hnorm{1}} + \norm{\dt v}{\Hnorm{1}} \norm{\sigma}{\Hnorm{2}} + \norm{v}{\Hnorm{2}}\norm{\dt \sigma}{\Hnorm{1}} + \norm{\dt \sigma}{\Hnorm{2}} \\
	& ~~~~ + \norm{v}{\Hnorm{3}} + \norm{v}{\Hnorm{2}} \norm{\sigma}{\Hnorm{3}} \bigr), \\
	& \norm{\dt w, \dz \dt w}{\Hnorm{1}} \leq C e^{2\norm{\sigma}{\Hnorm{2}}} \bigl( \norm{\sigma}{\Hnorm{3}}^3 + \norm{\dt \sigma}{\Hnorm{2}}^3 + 1 \bigr) \bigl( \norm{\dt^2 \sigma}{\Hnorm{1}}\\
	& ~~~~ + \norm{\dt v}{\Hnorm{2}} + \norm{\dt v}{\Hnorm{2}}\norm{\sigma}{\Hnorm{2}} + \norm{v}{\Hnorm{2}}\norm{\dt \sigma}{\Hnorm{2}} + \norm{\dt \sigma}{\Hnorm{2}}\\
	& ~~~~ + \norm{v}{\Hnorm{3}} + \norm{v}{\Hnorm{3}} \norm{\sigma}{\Hnorm{3}} \bigr),
	\end{aligned}
	\end{equation}
	for some generic positive constant $ C \in (0,\infty) $, independent of $ \varepsilon $. 
	In particular, 
	\begin{equation}\label{eqvl-functionals-t}
	\begin{gathered}
	\mathcal E_1 \leq E_1^2 \leq \mathfrak N(\mathcal E, \mathcal E_1),\\
	D_1 \leq \mathfrak N (\mathcal E, \mathcal E_1,1) \sqrt{\mathcal D_1} + \mathfrak N (\mathcal E,\mathcal E_1).
	\end{gathered}
	\end{equation}
\end{prop}
The proofs of Proposition \ref{lm:nonlinearity-t} and Proposition \ref{lm:vertical-velocity-t} are similar to those of Proposition \ref{lm:nonlinearity} and Proposition \ref{lm:vertical-velocity}, which are postponed to section \ref{subsec:proof-est-t}.

Now we apply the estimates in Proposition \ref{lm:nonlinearity-t}. Indeed, \eqref{est-101}, \eqref{est-1011}, \eqref{est-102}, \eqref{est-103}, \eqref{est-nonlinearity-i-t-1}, and \eqref{est-nonlinearity-j-t} yield
\begin{equation*}
\dfrac{d}{dt} \mathcal E_1 + \mathcal D_1 \leq \mathfrak N(E, E_1) ( 1 + D + D_1).
\end{equation*}
After substituting \eqref{eqvl-functionals} and \eqref{eqvl-functionals-t} on the right-hand side of the above inequality and applying the Young inequality in the resultant, one gets
\begin{equation}\label{ineq:total-1-2}
\dfrac{d}{dt} \mathcal E_1 + \mathcal D_1 \leq \mathfrak N( \mathcal E, \mathcal E_1,1)+ \dfrac{1}{2}\mathcal D_1 + \mathcal D .
\end{equation}
Thus, there is some constant $ T^{**} \in (0, T^*] $, independent of $ \varepsilon $, provided that $ \mathcal E(0), \mathcal E_{1}(0) $ are finite, such that
\begin{equation}\label{uniform-est-1-t}
\sup_{0\leq t \leq T^{**}} \mathcal E_1(t) + \int_0^{T^{**}} \mathcal D_1(t) \leq C < \infty,
\end{equation}
where $ C \in (0,\infty) $ depends on the initial data but is independent of $\varepsilon $. Here we have also used \eqref{uniform-est-1}. 
It follows from \eqref{uniform-est-1},  \eqref{eqvl-functionals-t}, and \eqref{uniform-est-1-t}, that
\begin{equation}\label{uniform-est-3}
\sup_{0\leq t \leq T^{**}} E_1^2(t) + \int_0^{T^{**}} D_1^2(t) \,dt < \infty, ~~~~ \text{for some} ~ T^{**} \in (0,T^*].
\end{equation}

\subsection{Convergence rates in terms of $ \varepsilon $}\label{subsec:converging-rate}

Our goal in this section is to estimate the convergence rates of $ \dlsigma, \dlv, \dlw $ to zero, as $ \varepsilon \rightarrow 0^+ $. 

In addition to \eqref{uniform-est-2} and \eqref{uniform-est-3}, let $ (\sigma_p, v_p, w_p) $ be a solution to the hydrostatic system \eqref{CPE}, satisfying
\begin{equation}\label{uniform-est-p}
\begin{gathered}
	\norm{\sigma_p}{L^\infty(0,T^{**};\Hnorm{4})} + \norm{\dt \sigma_p}{L^\infty(0,T^{**};\Hnorm{2})}+ \norm{v_p}{L^\infty(0,T^{**};\Hnorm{3})}  \\
	 + \norm{v_p}{L^2(0,T^{**};\Hnorm{4})} + \norm{w_p}{L^\infty(0,T^{**};\Hnorm{2})}  \leq C < \infty,
\end{gathered}
\end{equation}
for some constant $ C \in ( 0, \infty ) $. 
Such a solution can be obtained as the limit in Theorem \ref{thm:unifrom_hydrostatic}, restricted to the sub-interval of time $ \lbrack 0, T^{**} \rbrack \subset \lbrack 0, T^{*} \rbrack $, which was established in section \ref{sec:convergence-cpe}. Alternatively, one can refer to \cite{LT2018a} for the proof of existence of such a solution to the compressible primitive equations. 
Let $ \mathcal C \in (0,\infty) $, independent of $ \varepsilon $, be the bound given by \eqref{uniform-est-2}, \eqref{uniform-est-3}, and \eqref{uniform-est-p},
i.e., 
\begin{equation}\label{uniform-est-part-1}
\begin{aligned}
& \norm{v, \varepsilon w}{L^\infty(0, T^{**};\Hnorm{3} )} + \norm{\dt v, \dt (\varepsilon w)}{L^\infty(0, T^{**};\Hnorm{1} )} +
 \norm{v, \varepsilon w}{L^2(0, T^{**};\Hnorm{4} )} \\
&  + \norm{\dt v, \dt (\varepsilon w)}{L^2(0, T^{**};\Hnorm{2} )} + \norm{\dt^2 v, \dt^2(\varepsilon w)}{L^2(0, T^{**};\Lnorm{2})}  \\
& + \norm{\dt \sigma, \nablah \sigma, \dfrac{\dz \sigma}{\varepsilon}}{L^\infty(0,T^{**};H^2)} 
+ \norm{\sigma}{L^\infty(0,T^{**};\Hnorm{4})} \\
& + \norm{\dt^2 \sigma, \nablah \dt \sigma, \dfrac{\dz \dt \sigma}{\varepsilon}}{L^\infty(0,T^{**};L^2)} 
+ \norm{\dt \sigma}{L^\infty(0,T^{**};\Hnorm{2})} \\
& + \norm{\dt \sigma, \nablah \sigma, \dfrac{\dz \sigma}{\varepsilon}}{L^2(0,T^{**};H^3)} + \norm{\dt^2 \sigma, \nablah \dt \sigma, \dfrac{\dz \dt \sigma}{\varepsilon}}{L^2(0,T^{**};H^1)} \\
& + \norm{w,\dz w}{L^\infty(0,T^{**};H^2)} + \norm{\dt w,\dz \dt  w}{L^\infty(0,T^{**};L^2)} \\
& + \norm{w,\dz w}{L^2(0,T^{**};H^3)} + \norm{\dt w,\dz \dt  w}{L^2(0,T^{**};H^1)} \\
& + \norm{\sigma_p}{L^\infty(0,T^{**};\Hnorm{4})} + \norm{\dt \sigma_p}{L^\infty(0,T^{**};\Hnorm{2})} + \norm{v_p}{L^\infty(0,T^{**};\Hnorm{3})} \\
& + \norm{v_p}{L^2(0,T^{**};\Hnorm{4})} + \norm{w_p}{L^\infty(0,T^{**};\Hnorm{2})} < \mathcal C. 
\end{aligned}
\end{equation}
Using \eqref{uniform-est-part-1} and the triangle inequality, we have
\begin{equation}\label{uniform-est-part-2}
	\begin{gathered}
		\norm{\dlsigma}{L^\infty(0,T^{**};\Hnorm{4})} + \norm{\dt \dlsigma}{L^\infty(0,T^{**};\Hnorm{2})}+ \norm{\dlv}{L^\infty(0,T^{**};\Hnorm{3})}  \\
	 + \norm{\dlv}{L^2(0,T^{**};\Hnorm{4})} + \norm{\dlw}{L^\infty(0,T^{**};\Hnorm{2})}  < 2\mathcal C.
	\end{gathered}
\end{equation}
Noticing that $ \dz \sigma_p = 0 $ and therefore $ \dz \sigma = \dz \dlsigma $, 
\eqref{uniform-est-part-1} implies that
\begin{equation}\label{con-rate-sigma-dz}
\begin{aligned}
& \norm{\dz \dlsigma}{L^\infty(0,T^{**};H^2)} + \norm{\dz \dlsigma}{L^2(0,T^{**};H^3)} \\
& ~~~~ + \norm{\dz \dt \dlsigma}{L^\infty(0,T^{**}, L^2)}  + \norm{\dz\dt \dlsigma}{L^2(0,T^{**},H^1)}  \leq 2 \varepsilon \mathcal C = \mathcal O(\varepsilon).
\end{aligned}
\end{equation}
This gives us partial information of the convergence rates. 

Furthermore, we separate $ \dlsigma $ by its vertical average and fluctuation, defined as in \eqref{def:avg-flc}, i.e., 
\begin{equation}\label{def:dmpt-deltasigma}
	\dlsigma = \overline \dlsigma + \widetilde{\dlsigma}. 
\end{equation}
Then the following inequalities are ture:
\begin{align}
	\norm{\widetilde {\dlsigma}}{L^2} & \lesssim \norm{\partial_z \dlsigma}{L^2}, \label{ineq:fl-dlsigma-3}\\
	\norm{\nablah \widetilde {\dlsigma}}{L^2} & \lesssim \norm{\partial_z \nablah \dlsigma}{L^2}, \label{ineq:fl-dlsigma-1}\\
	\norm{\widetilde{\dt \dlsigma}}{L^2} & \lesssim \norm{\dz \dt \dlsigma}{L^2}. \label{ineq:fl-dlsigma-2}
\end{align}
The proof of \eqref{ineq:fl-dlsigma-3}--\eqref{ineq:fl-dlsigma-2} follows directly from applying the one-dimensional Poincar\'e inequality in the vertical direction, i.e., 
\begin{equation*}
	\norm{\widetilde{\phi}(\cdot)}{L^2(\mathbb T)} = \norm{{\phi}(\cdot) - \overline \phi }{L^2(\mathbb T)} \lesssim \norm{\partial_z \phi(\cdot)}{L^2(\mathbb T)}.
\end{equation*} 

We establish the rest of the required estimates in the following steps: the estimate of $ \overline{\dlsigma} $; the estimate of $ \dlv $; summary of the estimates; convergence rates via interpolation.

{\noindent\bf Step 1: }estimates of $ \overline{\dlsigma} $.
Taking the vertical-average of \subeqref{eq:difference}{1} leads to 
\begin{equation}\label{con-rate-001}
\dt \overline{\dlsigma} + \overline{v \cdot \nablah \dlsigma} + \overline{\dlv} \cdot \nablah \sigma_p + \dvh \overline{\dlv} + \overline{w \dz \dlsigma} = 0.
\end{equation}
Then we take the $ L^2 $-inner product of \eqref{con-rate-001} with $ 2 \overline{\dlsigma} $, it follows, 
\begin{equation}\label{con-rate-005}
\begin{aligned}
& \dfrac{d}{dt} \hnorm{\overline{\dlsigma}}{\Lnorm{2}}^2 = - 2 \int_{\Omega_h} \overline{v\cdot \nablah \dlsigma} \times \overline{\dlsigma} \idxh \\
& ~~~~ - 2 \int_{\Omega_h} (\overline{\dlv}\cdot\nablah \sigma_p + \dvh \overline{\dlv} + \overline{w\dz \dlsigma}) \times \overline{\dlsigma} \idxh =: N_1 + N_2.
\end{aligned}
\end{equation}
$ N_2 $ can be estimated, after applying the H\"older, Minkowski, and Sobolev embedding inequalities, as
\begin{align*}
& N_2 \lesssim \bigl(\norm{\dlv}{\Lnorm{2}} \hnorm{\nablah \sigma_p}{\Lnorm{\infty}} + \norm{\nabla \dlv}{\Lnorm{2}} + \norm{w}{\Lnorm{\infty}} \norm{\dz \dlsigma}{\Lnorm{2}} \bigr) \hnorm{\overline{\dlsigma}}{\Lnorm{2}} \\
& ~~~~ \lesssim \bigl(\norm{\dlv}{\Lnorm{2}} \hnorm{\nablah \sigma_p}{\Hnorm{2}} + \norm{\nabla \dlv}{\Lnorm{2}} + \norm{w}{\Hnorm{2}} \norm{\dz \dlsigma}{\Lnorm{2}} \bigr) \hnorm{\overline{\dlsigma}}{\Lnorm{2}} \\
& ~~~~ \lesssim \mathcal C \hnorm{\overline{\dlsigma}}{\Lnorm{2}} \norm{\dlv}{\Lnorm{2}} + \norm{\nabla\dlv}{\Lnorm{2}}\hnorm{\overline{\dlsigma}}{\Lnorm{2}} + \varepsilon \mathcal C \hnorm{\overline{\dlsigma}}{\Lnorm{2}},
\end{align*}
where we have substituted the bounds in \eqref{uniform-est-part-1} and \eqref{con-rate-sigma-dz}. 
On the other hand, $ N_1 $ can be written as, after substituting \eqref{def:dmpt-deltasigma} and integration by parts,
\begin{align*}
& N_1 = - 2 \int_{\Omega_h} \overline{v} \cdot \nablah \overline{\dlsigma} \times \overline{\dlsigma} \idxh 
- 2 \int_{\Omega_h} \overline{v \cdot \nablah \widetilde{\dlsigma}} \times \overline{\dlsigma} \idxh \\
& ~~~~  =  \int_{\Omega_h} \dvh \overline{v} \times  |\overline{\dlsigma}|^2\idxh 
- 2 \int_{\Omega_h} \overline{v \cdot \nablah \widetilde{\dlsigma}} \times \overline{\dlsigma} \idxh.
\end{align*}
Thus applying the H\"older, Minkowski, and Sobolev embedding inequalities and \eqref{ineq:fl-dlsigma-1} to the above identity yields
\begin{align*}
& N_1 \lesssim \norm{\nabla v}{\Lnorm{\infty}} \hnorm{\overline{\dlsigma}}{\Lnorm{2}}^2 + \norm{v}{\Lnorm{\infty}} \norm{\nablah \widetilde{\dlsigma}}{\Lnorm{2}} \hnorm{\overline{\dlsigma}}{\Lnorm{2}} \\
& ~~~~ \lesssim \norm{v}{\Hnorm{3}}\hnorm{\overline{\dlsigma}}{\Lnorm{2}}^2 + \norm{v}{\Hnorm{2}} \norm{\dz \dlsigma}{\Hnorm{1}} \hnorm{\overline{\dlsigma}}{\Lnorm{2}} \lesssim \mathcal C \hnorm{\overline{\dlsigma}}{\Lnorm{2}}^2 + \varepsilon \mathcal C \hnorm{\overline{\dlsigma}}{\Lnorm{2}},
\end{align*}
where in the last inequality we have substituted the bounds in \eqref{uniform-est-part-1} and \eqref{con-rate-sigma-dz}. 
Hence, \eqref{con-rate-005} implies that
\begin{equation}\label{con-rate-average-dlsigma}
\dfrac{d}{dt}\hnorm{\overline{\dlsigma}}{\Lnorm{2}}^2 \lesssim \mathcal C \bigl(\hnorm{\overline{\dlsigma}}{\Lnorm{2}}^2 + \norm{\dlv}{\Lnorm{2}}^2 \bigr) + \varepsilon \mathcal C \hnorm{\overline{\dlsigma}}{\Lnorm{2}} + \norm{\nabla \dlv}{\Lnorm{2}}\hnorm{\overline{\dlsigma}}{\Lnorm{2}}.
\end{equation}

{\noindent\bf Step 2: }estimates of $ \dlv $.

Taking the $ L^2 $-inner product of \subeqref{eq:difference}{2} with $ 2 \dlv $ leads to
\begin{equation}\label{con-rate-006}
\begin{aligned}
& \dfrac{d}{dt}\norm{\dlv}{\Lnorm{2}}^2 + 2 \norm{\nabla \dlv}{\Lnorm{2}}^2 = \int ( \dvh v + \dz w) |\dlv|^2 \idx\\
& ~~~~ - 2 \int (\dlv \cdot \nablah v_p) \cdot\dlv \idx + 2 \int \dlsigma ( \dvh \dlv ) \idx - \int \dlw ( \dz v_p  \cdot \dlv )\idx \\
& ~~~~=: N_3 + N_4 + N_5 + N_6. 
\end{aligned}
\end{equation}
Again, applying the H\"older and Sobolev embedding inequalities leads to
\begin{align*}
& N_3 + N_4 + N_5 \lesssim (\norm{\dvh v}{\Lnorm{\infty}}+\norm{\dz w}{\Lnorm{\infty}}+ \norm{\nablah v_p}{\Lnorm{\infty}}) \norm{\dlv}{\Lnorm{2}}^2 \\
& ~~~~ + \norm{\dlsigma}{\Lnorm{2}}\norm{\dvh \dlv}{\Lnorm{2}} \lesssim (\norm{v}{\Hnorm{3}}+\norm{\dz w}{\Hnorm{2}} + \norm{v_p}{\Hnorm{3}}) \norm{\dlv}{\Lnorm{2}}^2 \\
& ~~~~ + \norm{\dlsigma}{\Lnorm{2}}\norm{\nabla \dlv}{\Lnorm{2}} \lesssim \mathcal C \norm{\dlv}{\Lnorm{2}}^2 + (\hnorm{\overline\dlsigma}{\Lnorm{2}} + \varepsilon \mathcal C )\norm{\nabla \dlv}{\Lnorm{2}},
\end{align*}
where we have substituted \eqref{def:dmpt-deltasigma}, \eqref{ineq:fl-dlsigma-3} and the bounds in \eqref{uniform-est-part-1} and \eqref{con-rate-sigma-dz}. To estimate $ N_6 $, we first need to substitute \eqref{id:vertical-velocity-1} for $ \dlw $, which leads to
\begin{align*}
& N_6 = \int \biggl\lbrack e^{-\dlsigma} \biggl( \int_0^z \bigl\lbrack e^{\dlsigma} (\dt \dlsigma + v \cdot\nablah \dlsigma + \dlv \cdot \nablah \sigma_p + w_p\dz \dlsigma + \dvh \dlv ) \bigr\rbrack \,dz' \biggr) \\
& ~~~~ \times ( \dz v_p \cdot \dlv ) \biggr\rbrack \idx = \int \biggl\lbrack e^{-\dlsigma} \biggl(\int_0^z e^{\dlsigma} \dt \dlsigma \,dz'\biggr) \times ( \dz v_p \cdot \dlv )  \biggr\rbrack \idx \\
& ~~~~ + \int \biggl\lbrack e^{-\dlsigma} \biggl( \int_0^z e^{\dlsigma} (v\cdot \nablah \dlsigma) \,dz' \biggr) \times ( \dz v_p \cdot \dlv ) \biggr\rbrack \idx \\
& ~~~~ + \int \biggl\lbrack e^{-\dlsigma} \biggl(\int_0^z e^{\dlsigma} (\dlv \cdot \nablah \sigma_p + w_p \dz \dlsigma + \dvh \dlv) \,dz' \biggr) \times ( \dz v_p \cdot \dlv ) \biggr\rbrack \idx \\
& ~~~~ =: N_6' + N_6'' + N_6'''.
\end{align*}

$ N_6''' $ can be estimated, after applying the H\"older, Minkowski, and Sobolev embedding inequalities, as
\begin{align*}
& N_6''' \lesssim e^{2\norm{\dlsigma}{\Lnorm{\infty}}} \bigl( \norm{\dlv}{\Lnorm{2}} \norm{\nablah \sigma_p}{\Lnorm{\infty}}+ \norm{w_p}{\Lnorm{\infty}} \norm{\dz \dlsigma}{\Lnorm{2}} + \norm{\nablah \dlv}{\Lnorm{2}} \bigr) \\
& ~~~~ \times \norm{\dz v_p}{\Lnorm{\infty}} \norm{\dlv}{\Lnorm{2}} \lesssim e^{2\norm{\dlsigma}{\Hnorm{2}}} \bigl(\norm{\dlv}{\Lnorm{2}}\norm{\sigma_p}{\Hnorm{3}} + \norm{w_p}{\Hnorm{2}}\norm{\dz \dlsigma}{\Lnorm{2}}\\
& ~~~~ + \norm{\nabla \dlv}{\Lnorm{2}} \bigr) \times \norm{v_p}{\Hnorm{3}}\norm{\dlv}{\Lnorm{2}} \lesssim e^{2\mathcal C} \mathcal C \norm{\dlv}{\Lnorm{2}}^2 + \varepsilon e^{2\mathcal C} \mathcal C^2  \norm{\dlv}{\Lnorm{2}} \\
& ~~~~ + e^{2\mathcal C} \mathcal C \norm{\nabla \dlv}{\Lnorm{2}}\norm{\dlv}{\Lnorm{2}},
\end{align*}
where we have substituted the bounds in \eqref{uniform-est-part-1}, \eqref{uniform-est-part-2}, and \eqref{con-rate-sigma-dz}.
On the other hand, applying integration by parts in 
$ N_6'' $ leads to
\begin{align*}
& N_6'' = - \int \biggl\lbrack \biggl(\int_0^z (e^{\dlsigma} \dlsigma v ) \,dz'\biggr) \cdot \nablah (e^{-\dlsigma} \dz v_p \cdot \dlv ) \biggr\rbrack \idx \\
& ~~~~  - \int \biggl\lbrack \biggl(\int_0^z \dlsigma \dvh (e^{\dlsigma} v)  \,dz'\biggr) \times (e^{-\dlsigma} \dz v_p \cdot \dlv ) \biggr\rbrack\idx.
\end{align*}
Then, after expanding every term in the above expression and applying the H\"older, Minkowski, and Sobolev embedding inequalities, one can derive
\begin{align*}
& N_6'' \lesssim e^{2\norm{\dlsigma}{\Lnorm{\infty}}}  \norm{v}{\Lnorm{\infty}}\norm{\dz v_p}{\Lnorm{\infty}}\norm{\dlsigma}{\Lnorm{2}} \norm{\nabla\dlv}{\Lnorm{2}} \\
& ~~~~ + e^{2\norm{\dlsigma}{\Lnorm{\infty}}}  \bigl(  \norm{v}{\Lnorm{\infty}}\norm{\nablah\dz v_p}{\Lnorm{\infty}} + \norm{v}{\Lnorm{\infty}} \norm{\dz v_p}{\Lnorm{\infty}} \norm{\nablah \dlsigma}{\Lnorm{\infty}} \\
& ~~~~ + \norm{\nablah v}{\Lnorm{\infty}} \norm{\dz v_p}{\Lnorm{\infty}}  \bigr)\norm{\dlsigma}{\Lnorm{2}} \norm{\dlv}{\Lnorm{2}}\\
& ~~~~ \lesssim e^{4\mathcal C} (\mathcal C^3 + \mathcal C^2 + \norm{v_p}{\Hnorm{4}}^2) (\hnorm{\overline\dlsigma}{\Lnorm{2}} + \varepsilon \mathcal C ) ( \norm{\nabla\dlv}{\Lnorm{2}} + \norm{\dlv}{\Lnorm{2}} ),
\end{align*}
where in the last inequality we have substituted \eqref{def:dmpt-deltasigma}, \eqref{ineq:fl-dlsigma-3}, \eqref{uniform-est-part-1}, \eqref{uniform-est-part-2}, and \eqref{con-rate-sigma-dz}. 

In order to estimate $ N_6' $, we first rewrite $ \dt \dlsigma $ into its average and fluctuation parts; 
that is, 
\begin{equation}\label{con-rate-201}
\begin{aligned}
& \dt \dlsigma = \overline{\dt \dlsigma} + \widetilde{\dt \dlsigma} 
= - \bigl( \overline{v \cdot \nablah \dlsigma} + \overline{\dlv} \cdot \nablah \sigma_p \\
& ~~~~ ~~~~ + \dvh \overline{\dlv} 
+ \overline{w \dz \dlsigma} \bigr) 
+ \widetilde{\dt\dlsigma},
\end{aligned}\end{equation}
where we have substituted \eqref{con-rate-001}. 
Thus $ N_6' $ can be written as
\begin{align*}
& N_6' = - \int \biggl\lbrack e^{-\dlsigma} \biggl( \int_0^z e^{\dlsigma} \overline{v \cdot \nablah \dlsigma} \,dz' \biggr) \times (\dz v_p \cdot \dlv)  \biggr\rbrack \idx \\
& ~~~~ - \int \biggl\lbrack e^{-\dlsigma} \biggl( \int_0^z e^{\dlsigma} ( \overline{\dlv}\cdot \nablah \sigma_p + \dvh \overline{\dlv} + \overline{w \dz \dlsigma } ) \,dz' \biggr) \times (\dz v_p \cdot \dlv)  \biggr\rbrack \idx \\
& ~~~~ + \int \biggl\lbrack e^{-\dlsigma} \biggl( \int_0^z e^{\dlsigma} \widetilde{\dt\dlsigma} \,dz' \biggr) \times (\dz v_p \cdot \dlv)  \biggr\rbrack \idx =: N_{6,1}' + N_{6,2}' + N_{6,3}'.
\end{align*} 
$ N_{6,1}' $ and $ N_{6,2}' $ can be estimated in the same way as that of $ N_6'' $ and $ N_6''' $, respectively. $ N_{6,3}' $ can be estimated, after applying the H\"older, Minkowski, Sobolev embedding inequalities and \eqref{ineq:fl-dlsigma-2}, as
\begin{align*}
& N_{6,3}' \lesssim e^{2\norm{\dlsigma}{\Lnorm{\infty}}} \norm{\dz v_p}{\Lnorm{\infty}} \norm{\dz \dt \dlsigma}{\Lnorm{2}} \norm{\dlv }{\Lnorm{2}} \\
& ~~~~ \lesssim e^{2\norm{\dlsigma}{\Hnorm{2}}} \norm{ v_p}{\Hnorm{3}} \norm{\dz \dt \dlsigma}{\Lnorm{2}} \norm{\dlv }{\Lnorm{2}} \lesssim \varepsilon e^{4\mathcal C} \mathcal C^2 \norm{\dlv}{\Lnorm{2}},
\end{align*}
where we have substituted the bounds in \eqref{uniform-est-part-1}, \eqref{uniform-est-part-2}, and \eqref{con-rate-sigma-dz} in the last inequality. Consequently, after summing up the above estimates, \eqref{con-rate-006} implies that
\begin{equation}\label{con-rate-dlv}
\begin{aligned}
& \dfrac{d}{dt}\norm{\dlv}{\Lnorm{2}}^2 + 2 \norm{\nabla \dlv}{\Lnorm{2}}^2 \lesssim e^{4\mathcal C} ( 1 +   \mathcal C^3 + \norm{v_p}{\Hnorm{4}}^2 ) ( \norm{\dlv}{\Lnorm{2}} + \hnorm{\overline{\dlsigma}}{\Lnorm{2}} + \varepsilon) \\
& ~~~~ ~~~~ ~~~~ \times ( \norm{\dlv}{\Lnorm{2}} + \norm{\nabla\dlv}{\Lnorm{2}}).
\end{aligned}
\end{equation}

{\noindent\bf Step 3: }summary of the estimates and convergence rates.
\eqref{con-rate-average-dlsigma} and \eqref{con-rate-dlv} imply, after applying the Young inequality, 
\begin{equation}\label{con-rate-dlsigma-dlv}
\dfrac{d}{dt} \bigl( \hnorm{\overline\dlsigma}{\Lnorm{2}}^2 + \norm{\dlv}{\Lnorm{2}}^2 \bigr) + \norm{\nabla\dlv}{\Lnorm{2}}^2 \leq C e^{4\mathcal C} (1+\mathcal C^3+ \norm{v_p}{\Hnorm{4}}^2) ( \hnorm{\overline{\dlsigma}}{\Lnorm{2}}^2 + \norm{\dlv}{\Lnorm{2}}^2 + \varepsilon^2  ).
\end{equation}
Therefore, applying the Gr\"onwall inequality to \eqref{con-rate-dlsigma-dlv} yields
\begin{equation}\label{con-rate-dlsigma-dlv-2}
\begin{gathered}
\sup_{0\leq t\leq T^{**}} \bigl( \hnorm{\overline\dlsigma(t)}{\Lnorm{2}}^2 + \norm{\dlv(t)}{\Lnorm{2}}^2 \bigr) + \int_0^{T^{**}} \norm{\nabla\dlv(t)}{\Lnorm{2}}^2 \,dt \\
\leq \mathfrak N (\mathcal C, T^{**}) ( \hnorm{\overline\dlsigma_0}{\Lnorm{2}}^2 + \norm{\dlv_0}{\Lnorm{2}}^2 + \varepsilon^2 ),
\end{gathered}
\end{equation}
where we have used the fact $ \norm{v_p}{L^2(0,T^{**};H^4)} < \mathcal C $. 
Together with \eqref{con-rate-sigma-dz} and \eqref{ineq:fl-dlsigma-3}, \eqref{con-rate-dlsigma-dlv-2} implies 
\begin{equation}\label{con-rate-dlsigma-dlv-3}
\begin{gathered}
\norm{\dlsigma}{L^\infty(0,T^{**};L^2)} + \norm{\dlv}{L^\infty(0,T^{**};L^2)} + \norm{\dlv}{L^2(0,T^{**};H^1)} \\
~~~~ \leq \mathfrak N (\mathcal C,T^{**}) ( \norm{\dlsigma_0}{\Lnorm{2}} + \norm{\dlv_0}{\Lnorm{2}} + \varepsilon ) = \mathcal O(\varepsilon),
\end{gathered}
\end{equation}
provided that $\norm{\dlsigma_0}{\Lnorm{2}} + \norm{\dlv_0}{\Lnorm{2}} = \mathcal O (\varepsilon) $, which is assumption \eqref{WP} for well-prepared initial data.

{\noindent\bf Step 4: }convergence rates via interpolation. 
Notice that, the Gagliardo-Nirenberg interpolation inequality implies
\begin{equation*}
\norm{\dlsigma}{\Hnorm{1}} \lesssim \norm{\dlsigma}{\Lnorm{2}}^{3/4} \norm{\dlsigma}{\Hnorm{4}}^{1/4}, ~~ \norm{\dlv}{\Hnorm{1}} \lesssim \norm{\dlv}{\Lnorm{2}}^{2/3} \norm{\dlv}{\Hnorm{3}}^{1/3}.
\end{equation*}
Therefore, together with \eqref{uniform-est-part-2} and \eqref{con-rate-dlsigma-dlv-3}, it follows,
\begin{equation}\label{con-rate-dlsigma-1}
\begin{gathered}
\norm{\dlsigma}{L^\infty(0,T^{**},H^1)} \lesssim \mathfrak N (\mathcal{C},T^{**}) \varepsilon^{3/4} = \mathcal O(\varepsilon^{3/4}), \\
\norm{\dlv}{L^\infty(0,T^{**}, H^1)} \lesssim \mathfrak N (\mathcal C,T^{**}) \varepsilon^{2/3} = \mathcal O (\varepsilon^{2/3}).
\end{gathered}\end{equation}

From \eqref{con-rate-201}, one can derive, after applying the H\"older, Minkowski, and Sobolev embedding inequalities,
\begin{equation}\label{con-rate-202}
\begin{aligned}
& \norm{\dt \dlsigma}{\Lnorm{2}} \lesssim \norm{v}{\Lnorm{\infty}} \norm{\nablah \dlsigma}{\Lnorm{2}} + \norm{\dlv}{\Lnorm{2}} \norm{\nablah \sigma_p}{\Lnorm{\infty}} + \norm{\nablah \dlv}{\Lnorm{2}} \\
& ~~~~ + \norm{w}{\Lnorm{\infty}} \norm{\dz \dlsigma}{\Lnorm{2}} + \norm{\dz \dt \dlsigma }{\Lnorm{2}} \lesssim \norm{v}{\Hnorm{2}} \norm{\dlsigma}{\Hnorm{1}} + \norm{\sigma_p}{\Hnorm{3}} \norm{\dlv}{\Lnorm{2}}\\
& ~~~~ + \norm{\dlv}{\Hnorm{1}} + \norm{w}{\Hnorm{2}}\norm{\dz\dlsigma}{\Lnorm{2}} + \norm{\dz \dt \dlsigma}{\Lnorm{2}}.
\end{aligned}\end{equation}
Thus \eqref{uniform-est-part-1}, \eqref{con-rate-sigma-dz}, \eqref{con-rate-dlsigma-dlv-3} \eqref{con-rate-dlsigma-1}, and \eqref{con-rate-202} yield,
\begin{equation}\label{con-rate-dlsigma-2}
\begin{gathered}
\norm{\dt \dlsigma}{L^\infty(0,T^{**},L^2)} \lesssim \mathfrak N(\mathcal C,T^{**}) \varepsilon^{2/3} = \mathcal O(\varepsilon^{2/3}), \\ 
\norm{\dt \dlsigma}{L^2(0,T^{**},L^2)} \lesssim \mathfrak N(\mathcal C,T^{**}) \varepsilon^{3/4} = \mathcal O(\varepsilon^{3/4}).
\end{gathered}\end{equation}
Similarly, from \eqref{id:vertical-velocity-1}, after applying the H\"older, Minkowski, and Sobolev embedding inequalities, one can derive,
\begin{equation}\label{con-rate-203}
\begin{aligned}
& \norm{\dlw}{\Lnorm{2}} \lesssim e^{2\norm{\dlsigma}{\Hnorm{2}}} \bigl( \norm{\dt \dlsigma}{\Lnorm{2}} +  \norm{v}{\Hnorm{2}} \norm{\dlsigma}{\Hnorm{1}} + \norm{\sigma_p}{\Hnorm{3}} \norm{\dlv}{\Lnorm{2}}\\
& ~~~~ ~~~~ + \norm{\dlv}{\Hnorm{1}} + \norm{w_p}{\Hnorm{2}}\norm{\dz\dlsigma}{\Lnorm{2}} \bigr).
\end{aligned}
\end{equation}
Therefore, substituting \eqref{uniform-est-part-1}, \eqref{uniform-est-part-2}, \eqref{con-rate-sigma-dz}, \eqref{con-rate-dlsigma-dlv-3}, \eqref{con-rate-dlsigma-1}, and \eqref{con-rate-dlsigma-2} to \eqref{con-rate-203} yields,
\begin{equation}\label{con-rate-204}
\begin{gathered}
\norm{\dlw}{L^\infty(0,T^{**};L^2)} \lesssim \mathfrak N (\mathcal C,T^{**}) \varepsilon^{2/3} = \mathcal O(\varepsilon^{2/3}) ,\\
\norm{\dlw}{L^2(0,T^{**};L^2)}\lesssim \mathfrak N (\mathcal C,T^{**}) \varepsilon^{3/4} = \mathcal O(\varepsilon^{3/4}) .
\end{gathered}
\end{equation}

\subsection{Proofs of Proposition \ref{lm:nonlinearity-t} and Proposition \ref{lm:vertical-velocity-t}}\label{subsec:proof-est-t}

\begin{proof}[Proof of Proposition \ref{lm:nonlinearity-t}]
	The proof is similar to that of Proposition \ref{lm:nonlinearity}. We list the estimates for readers' convenience:
	\begin{align*}
	& \norm{\dt \mathcal I_1}{\Lnorm{2}} \lesssim \norm{\dt v}{\Lnorm{6}}\norm{\nabla v, \nabla(\varepsilon w)}{\Lnorm{3}} + \norm{v}{\Lnorm{\infty}} \norm{\nabla \dt v, \nabla\dt(\varepsilon w)}{\Lnorm{2}} \\
	& ~~~~ ~~~~ \lesssim \norm{\dt v, \dt(\varepsilon w)}{\Hnorm{1}}\norm{v, \varepsilon w}{\Hnorm{2}}\lesssim E_1 E, \\
	& \norm{\dt \mathcal I_2}{\Lnorm{2}} \lesssim \norm{\dt w}{\Lnorm{2}}\norm{\dz(\varepsilon w),\dz v}{\Lnorm{\infty}} + \norm{w}{\Lnorm{\infty}}\norm{\dz \dt (\varepsilon w), \dz \dt v}{\Lnorm{2}} \\
	& ~~~~ ~~~~ \lesssim \norm{\dt w}{\Lnorm{2}} \norm{\varepsilon w, v}{\Hnorm{3}} + \norm{w }{\Hnorm{2}} \norm{\dt(\varepsilon w), \dt v}{\Hnorm{1}} \lesssim E_1 E, 
	\end{align*}
	where we have applied the H\"older and Sobolev embedding inequalities. Similarly, we have,
	\begin{align*}
	& \norm{\dt \mathcal J_1}{\Lnorm{2}} \lesssim \norm{\dt^2 v, \dt \Delta v }{\Lnorm{2}}\norm{\nablah \sigma}{\Lnorm{\infty}} + \norm{\dt v, \Delta v}{\Lnorm{3}} \norm{\nablah \dt \sigma}{\Lnorm{6}} \\
	& ~~~~ + \norm{\dt \nabla v}{\Lnorm{2}}\norm{\nablah \nabla \sigma}{\Lnorm{\infty}} + \norm{\nabla v}{\Lnorm{\infty}} \norm{\nablah\nabla\dt \sigma}{\Lnorm{2}}\\
	& ~~~~  \lesssim E (E + E_1 + D + D_1),   \\
	& \norm{\dt \mathcal J_2}{\Lnorm{2}} \lesssim \norm{\dt^2(\varepsilon w), \dt \Delta (\varepsilon w)}{\Lnorm{2}}\norm{\dfrac{\dz\sigma}{\varepsilon}}{\Lnorm{\infty}} + \norm{\dt(\varepsilon w), \Delta(\varepsilon w)}{\Lnorm{3}}\norm{\dfrac{\dz \dt \sigma}{\varepsilon}}{\Lnorm{6}} \\
	& ~~~~  + \norm{\nabla\dt (\varepsilon w)}{\Lnorm{3}} \norm{\dfrac{\nabla\dz \sigma}{\varepsilon}}{\Lnorm{6}} + \norm{\nabla(\varepsilon w)}{\Lnorm{2}}\norm{\dfrac{\nabla \dz \dt \sigma}{\varepsilon}}{\Lnorm{2}} \\
	& ~~~~  \lesssim E (E + E_1 + D + D_1 ), \\
	& \norm{\dt \mathcal J_3}{\Lnorm{2}} \lesssim \norm{\nabla \dt v}{\Lnorm{2}} \norm{\nablah v}{\Lnorm{\infty}} 
	  + \norm{\nabla v}{\Lnorm{\infty}} \norm{\nablah \dt v}{\Lnorm{2}} 
	+ \norm{\dt v}{\Lnorm{\infty }}\norm{\nabla\nablah v}{\Lnorm{2}}
 	\\ & ~~~~ + \norm{v}{\Lnorm{\infty}} \norm{\nabla\nablah \dt v}{\Lnorm{2}}  + \norm{\nabla\dt w}{\Lnorm{2}} \norm{\dz v}{\Lnorm{\infty}} 
	 + \norm{\nabla w}{\Lnorm{\infty}} \norm{\dz \dt v}{\Lnorm{2}} 	\\ & ~~~~ + \norm{\dt w}{\Lnorm{3}} \norm{\nabla\dz v}{\Lnorm{6}}  + \norm{w}{\Lnorm{\infty}} \norm{\nabla \dz \dt v}{\Lnorm{2}} \\
	& ~~~~  \lesssim (E + E_1 ) ( E + E_1 + D + D_1),  \\
	& \norm{\dt \mathcal J_4}{\Lnorm{2}} \lesssim \norm{\dz \dt v}{\Lnorm{3}} \norm{\nablah w}{\Lnorm{6}} + \norm{\dz v}{\Lnorm{\infty}}\norm{\nablah \dt w}{\Lnorm{2}} 
	 + \norm{ \dt v }{\Lnorm{3}}\norm{\nablah \dz w}{\Lnorm{6}} 
	 \\ & ~~~~ + \norm{v}{\Lnorm{\infty}} \norm{\nablah \dz \dt w}{\Lnorm{2}} 
	 + \norm{\dz\dt w}{\Lnorm{2}}\norm{\dz w}{\Lnorm{\infty}} + \norm{\dz w}{\Lnorm{\infty}}\norm{\dz \dt w, }{\Lnorm{2}} \\
	& ~~~~ + \norm{\dt w }{\Lnorm{3}} \norm{\dz^2 w}{\Lnorm{6}} + \norm{w}{\Lnorm{\infty}}\norm{\dz^2 \dt w}{\Lnorm{2}}\\
	& ~~~~  \lesssim (E + E_1 ) ( E + E_1 + D + D_1).
	\end{align*}
	This completes the proof.
\end{proof}

\begin{proof}[Proof of Proposition \ref{lm:vertical-velocity-t}]
	We follow similar steps as in the proof of Proposition \ref{lm:vertical-velocity}. Recalling $ \Xi $ in \eqref{id:vertical-velocity-notation}, applying a temporal derivative to \eqref{id:vertical-velocity-0} leads to
	\begin{equation*}
	\dt w = - e^{-\sigma} \int_0^z e^\sigma (\dt \sigma \Xi + \dt \Xi) \,dz'+ e^{-\sigma} \dt \sigma  \int_0^z e^\sigma \Xi \,dz'.
	\end{equation*}
	Similarly, 
	\begin{align*}
	& \dz \dt w = e^{-\sigma} \dz \sigma \int_0^z e^\sigma(\dt\sigma \Xi + \dt \Xi)\,dz' + e^{-\sigma}(\dz \dt \sigma - \dz \sigma \dt \sigma ) \int_0^z e^\sigma \Xi \,dz' \\
	& ~~~~ ~~~~ -\dt\Xi, \\
	& \ddh \dt w = e^{-\sigma} \ddh \sigma \int_0^z e^\sigma (\dt \sigma \Xi + \dt \Xi ) \,dz' + e^{-\sigma}(\ddh \dt \sigma - \ddh \sigma \dt \sigma )\int_0^z e^\sigma \Xi \,dz'\\
	& ~~~~ ~~~~ - e^{-\sigma} \int_0^z e^\sigma (  \ddh \sigma \dt \sigma \Xi + \ddh \sigma \dt \Xi + \ddh \dt \sigma \Xi + \dt \sigma \ddh \Xi + \ddh \dt \Xi  )\,dz'\\
	& ~~~~ ~~~~ + e^{-\sigma} \dt \sigma \int_0^z e^\sigma (\ddh \sigma \Xi + \ddh \Xi ) \,dz' , \\
	& \partial_{hz} \dt w = e^{-\sigma}\dz \sigma \int_0^z e^{\sigma} ( \ddh \sigma \dt \sigma \Xi + \ddh \sigma \dt \Xi + \dt \sigma \ddh \Xi + \ddh \dt \sigma \Xi + \ddh \dt \Xi) \,dz' \\
	& ~~~~ ~~~~ + e^{-\sigma}(\partial_{hz} \sigma - \ddh \sigma \dz \sigma)\int_0^z e^\sigma (\dt \sigma \Xi + \dt \Xi)\,dz' + e^{-\sigma} (\dz \dt \sigma - \dz \sigma \dt \sigma) \\
	& ~~~~ ~~~~ ~~ \times \int_0^z e^\sigma( \ddh \sigma \Xi + \ddh \Xi )\,dz' + e^{-\sigma}( \partial_{hz}\dt \sigma - \partial_{hz} \sigma \dt \sigma - \dz \sigma \ddh \dt \sigma \\
	& ~~~~ ~~~~ ~~ - \ddh \sigma \dz \dt \sigma + \ddh \sigma \dz \sigma \dt \sigma) \int_0^z e^\sigma \Xi \,dz' - \ddh \dt \Xi.
	\end{align*}
	
	Therefore, after applying the H\"older, Minkowski, and Sobolev embedding inequalities, one obtains,
	\begin{align}
	& \label{est-w-101}
	\begin{aligned}
	& \norm{\dt w, \dz \dt w}{\Lnorm{2}} \lesssim e^{2\norm{\sigma}{\Hnorm{2}}} \bigl(\norm{\sigma}{\Hnorm{3}} + 1 \bigr) \bigl(\norm{\dt^2 \sigma}{\Lnorm{2}} + \norm{\dt v}{\Hnorm{1}} \\
	& ~~~~ ~~~~ + \norm{\dt v}{\Hnorm{1}} \norm{\sigma}{\Hnorm{2}} + \norm{v}{\Hnorm{2}} \norm{\dt \sigma }{\Hnorm{1}} + \norm{\dt \sigma}{\Hnorm{1}} (\norm{\dt \sigma}{\Hnorm{2}} \\
	& ~~~~ ~~~~ + \norm{v}{\Hnorm{3}} + \norm{v}{\Hnorm{2}}\norm{\sigma}{\Hnorm{3}} )  \bigr),
	\end{aligned}\\
	& \label{est-w-102}
	\begin{aligned}
	& \norm{\ddh \dt w, \partial_{hz} \dt w }{\Lnorm{2}} \lesssim e^{2\norm{\sigma}{\Hnorm{2}}} \bigl( \norm{\sigma}{\Hnorm{3}}^2 + 1 \bigr) \bigl( \norm{\dt^2 \sigma}{\Hnorm{1}} + \norm{\dt v}{\Hnorm{2}} \\
	& ~~~~ ~~~~ + \norm{\dt v}{\Hnorm{2}} \norm{\sigma}{\Hnorm{2}} + \norm{v}{\Hnorm{2}}\norm{\dt \sigma}{\Hnorm{2}} + ( \norm{\dt \sigma}{\Hnorm{2}} + 1) \\
	& ~~~~ ~~~~ \times ( \norm{\dt \sigma }{\Hnorm{2}} + \norm{v}{\Hnorm{3}} + \norm{v}{\Hnorm{3}}\norm{\sigma}{\Hnorm{3}} )\bigr).
	\end{aligned}
	\end{align}
	This completes the proof of Proposition \ref{lm:vertical-velocity-t}. 
\end{proof}

\section*{Acknowledgements}

This work was supported in part by the Einstein Stiftung/Foundation - Berlin, through the Einstein Visiting Fellow Program. 
X.L. and E.S.T. would like to thank the Isaac Newton Institute for Mathematical Sciences, Cambridge, for support and hospitality during the programme {\it``Mathematical aspects of turbulence: where do we stand?"}, where part of the work on this paper was undertaken.
This work was supported in part by EPSRC grant no EP/R014604/1.
X.L.'s work was partially supported by a grant from the Simons Foundation, during his visit to the Isaac Newton Institute for Mathematical Sciences. The research of E.S.T. has benefited from the inspiring environment of the CRC 1114 ``Scaling Cascades in Complex Systems'', Project Number 235221301, Project C06, funded by Deutsche Forschungsgemeinschaft (DFG).

\bibliographystyle{plain}

\begin{thebibliography}{10}
	
	\bibitem{Azerad2001}
	Pascal Az{\'{e}}rad and Francisco Guill{\'{e}}n.
	\newblock {Mathematical justification of the hydrostatic approximation in the
		primitive equations of geophysical fluid dynamics}.
	\newblock {\em SIAM J. Math. Anal.}, 33(4):847--859, 2001.
	
	\bibitem{Bella2014}
	Peter Bella, Eduard Feireisl, and Anton{\'{i}}n Novotn{\'{y}}.
	\newblock {Dimension reduction for compressible viscous fluids}.
	\newblock {\em Acta Appl. Math.}, 134(1):111--121, 2014.
	
	\bibitem{Bresch2006}
	Didier Bresch and Beno\^{i}t Desjardins.
	\newblock {Stabilit\'{e} de solutions faibles globales pour les \'equations de
		Navier--Stokes compressible avec temp\'erature}.
	\newblock {\em C. R. Acad. Sci. Paris, Ser. I}, 343(3):219--224, 2006.
	
	\bibitem{Bresch2018}
	Didier Bresch and Pierre-Emmanuel Jabin.
	\newblock {Global existence of weak solutions for compressible Navier--Stokes
		equations: Thermodynamically unstable pressure and anisotropic viscous stress
		tensor}.
	\newblock {\em Ann. Math.}, 188(2):577--684, 2018.
	
	\bibitem{Bresch2019}
	Didier Bresch, Alexis Vasseur, and Cheng Yu.
	\newblock {Global existence of entropy-weak solutions to the compressible
		Navier--Stokes equations with non-linear density dependent viscosities}.
	\newblock {\em J. Eur. Math. Soc.}, 24(5):1791--1837, 2022.
	
	\bibitem{Cao2014a}
	Chongsheng Cao, Jinkai Li, and Edriss~S. Titi.
	\newblock {Global well-posedness of strong solutions to the 3D primitive
		equations with horizontal eddy diffusivity}.
	\newblock {\em J. Differential Equations}, 257(11):4108--4132, 2014.
	
	\bibitem{Cao2014b}
	Chongsheng Cao, Jinkai Li, and Edriss~S Titi.
	\newblock {Local and global well-posedness of strong solutions to the 3D
		primitive equations with vertical eddy diffusivity}.
	\newblock {\em Arch. Rational Mech. Anal.}, 214(1):35--76, 2014.
	
	\bibitem{Cao2016}
	Chongsheng Cao, Jinkai Li, and Edriss~S. Titi.
	\newblock {Global well-posedness of the three-dimensional primitive equations
		with only horizontal viscosity and diffusion}.
	\newblock {\em Communications on Pure and Applied Mathematics}, 69(8):1492--1531, 2016.
	
	\bibitem{Cao2017}
	Chongsheng Cao, Jinkai Li, and Edriss~S. Titi.
	\newblock {Global well-posedness of the 3D primitive equations with horizontal
		viscosity and vertical diffusivity}.
	\newblock {\em Physica D: Nonlinear Phenomena}, 412:132606, 2020. 
	
	\bibitem{Cao2016a}
	Chongsheng Cao, Jinkai Li, and Edriss~S. Titi.
	\newblock {Strong solutions to the 3D primitive equations with only horizontal
		dissipation: Near $ H^1 $ initial data}.
	\newblock {\em Journal of Functional Analysis}, 272(11):4606--4641, 2017.
	
	\bibitem{Cao2007}
	Chongsheng Cao and Edriss Titi.
	\newblock {Global well-posedness of the three-dimensional viscous primitive
		equations of large scale ocean and atmosphere dynamics}.
	\newblock {\em Ann. Math.}, 166(1):245--267, 2007.
	
	\bibitem{Cao2012}
	Chongsheng Cao and Edriss~S. Titi.
	\newblock {Global well-posedness of the 3D primitive equations with partial
		vertical turbulence mixing heat diffusion}.
	\newblock {\em Commun. Math. Phys.}, 310(2):537--568, 2012.
	
	\bibitem{Cho2004}
	Yonggeun Cho, Hi~Jun Choe, and Hyunseok Kim.
	\newblock {Unique solvability of the initial boundary value problems for
		compressible viscous fluids}.
	\newblock {\em J. Math. Pures Appl.}, 83(2):243--275, 2004.
	
	\bibitem{Cho2006a}
	Yonggeun Cho and Hyunseok Kim.
	\newblock {Existence results for viscous polytropic fluids with vacuum}.
	\newblock {\em J. Differential Equations}, 228(2):377--411, 2006.
	
	\bibitem{Cho2006c}
	Yonggeun Cho and Hyunseok Kim.
	\newblock {On classical solutions of the compressible Navier--Stokes equations
		with nonnegative initial densities}.
	\newblock {\em manuscripta Math.}, 120(1):91--129, 2006.
	
	\bibitem{Ersoy2012}
	Mehmet Ersoy and Timack Ngom.
	\newblock {Existence of a global weak solution to compressible primitive
		equations}.
	\newblock {\em C. R. Acad. Sci. Paris, Ser. I}, 350(7--8):379--382, 2012.
	
	\bibitem{Ersoy2011a}
	Mehmet Ersoy, Timack Ngom, and Mamadou Sy.
	\newblock {Compressible primitive equations: Formal derivation and stability of
		weak solutions}.
	\newblock {\em Nonlinearity}, 24(1):79--96, 2011.
	
	
	\bibitem{Feireisl2004}
	Eduard Feireisl.
	\newblock {\em {Dynamics of Viscous Compressible Fluids}}.
	\newblock Oxford Lecture Series in Mathematics and its Applications, 26. Oxford
	University Press, 2004.
	
	\bibitem{gao2022hydrostatic}
	Hongjun Gao, {\v{S}}{\'a}rka Ne{\v{c}}asov{\'a}, and Tong Tang.
	\newblock {On the hydrostatic approximation of compressible anisotropic Navier--Stokes equations--rigorous justification},
	\newblock {\rm Journal of Mathematical Fluid Mechanics},
  	24(3):86, {2022}.
	
	\bibitem{Gatapov2005}
	B.~V. Gatapov and Aleksandr Vasilevic Kazhikhov.
	\newblock {Existence of a global solution to one model problem of atmosphere
		dynamics}.
	\newblock {\em Siberian Mathematical Journal}, 46(5):805--812, 2005.
	
	\bibitem{GerardVaret2018}
	David Gerard-Varet, Nader Masmoudi, and Vlad Vicol.
	\newblock {Well-posedness of the hydrostatic Navier--Stokes equations}.
	\newblock {\em Anal. PDE}, 13(5):1417--1455, 2020.	
	
	\bibitem{GuillenGonzalez2001}
	Francisco~Guill\'{e}n-Gonz\'{a}lez, Nader~Masmoudi, and Mar\'ia \'Angeles Rodr\'{i}guez-Bellido.
	\newblock {Anisotropic estimates and strong solutions of the primitive
		equations}.
	\newblock {\em Differential and Integral Equations},
	14(11):1381--1408, 2001.
	
	\bibitem{Hoff1992}
	David Hoff.
	\newblock {Global well-posedness of the Cauchy problem for the Navier-Stokes
		equations of nonisentropic flow with discontinuous initial data}.
	\newblock {\em J. Differential Equations}, 95(1):33--74, 1992.
	
	\bibitem{hoff1991}
	David Hoff and Denis Serre.
	\newblock {The failure of continuous dependence on initial data for the
		Navier--Stokes equations of compressible flow}.
	\newblock {\em SIAM J. Appl. Math.}, 51(4):887--898, 1991.
	
	\bibitem{HuTemamZiane2003}
	Changbing Hu, Roger Temam, and Mohammed Ziane.
	\newblock {The primitive equations on the large scale ocean under the small
		depth hypothesis}.
	\newblock {\em Discrete and Continuous Dynamical Systems}, 9(1):97--131, 2003.
	
	\bibitem{Huang2016}
	Xiangdi Huang and Jing Li.
	\newblock {Existence and blowup behavior of global strong solutions to the two-dimensional barotropic compressible Navier--Stokes system with vacuum and large initial data}.
	\newblock {\em J. Math. Pures Appl.}, 106(1):123--154, 2016.
	
	\bibitem{Huang2018a}
	Xiangdi Huang and Jing Li.
	\newblock {Global classical and weak solutions to the three-dimensional full
	compressible Navier--Stokes system with vacuum and large oscillations}.
	\newblock {\em Arch. Rational Mech. Anal.}, 227(3):995--1059, 2018.
	
	\bibitem{HuangLiXin2012}
	Xiangdi Huang, Jing Li, and Zhouping Xin.
	\newblock {Global well-posedness of classical solutions with large oscillations
		and vacuum to the three-dimensional isentropic compressible Navier--Stokes
		equations}.
	\newblock {\em Commun. Pure Appl. Math.}, 65(4):549--585, 2012.
	
	\bibitem{Itaya1971}
	Nobutoshi Itaya.
	\newblock {On the Cauchy problems for the system of fundamental equations
		describing the movement of compressible viscous fluid}.
	\newblock {\em Kodai Math. Sem. Rep.}, 23(1):60--120, 1971.
	
	\bibitem{Jiu2018}
	Quansen Jiu, Mingjie Li, and Fengchao Wang.
	\newblock {Uniqueness of the global weak solutions to 2D compressible primitive
		equations}.
	\newblock {\em J. Math. Anal. Appl.}, 461(2):1653--1671, 2018.
	
	\bibitem{Choe2003}
	Hi~{Jun Choe} and Hyunseok Kim.
	\newblock {Strong solutions of the Navier--Stokes equations for isentropic
		compressible fluids}.
	\newblock {\em J. Differential Equations}, 190(2):504--523, 2003.
	
	\bibitem{Kukavica2007a}
	Igor Kukavica and Mohammed Ziane.
	\newblock {The regularity of solutions of the primitive equations of the ocean
		in space dimension three}.
	\newblock {\em C. R. Acad. Sci. Paris, Ser. I}, 345(5):257--260, 2007.
	
	\bibitem{Li2015a}
	Jing Li and Zhouping Xin.
	\newblock {Global existence of weak solutions to the barotropic compressible
		Navier--Stokes flows with degenerate viscosities}.
	\newblock Available at \href{https://arxiv.org/abs/1504.06826}{arXiv:1504.06826}, 2015.
	
	\bibitem{Li2017a}
	Jinkai Li and Edriss~S Titi.
	\newblock {Existence and uniqueness of weak solutions to viscous primitive
		equations for a certain class of discontinuous initial data}.
	\newblock {\em SIAM J. Math. Anal.}, 49(1):1--28, 2017.
	
	\bibitem{Li2017}
	Jinkai Li and Edriss~S. Titi.
	\newblock {The primitive equations as the small aspect ratio limit of the
		Navier--Stokes equations: Rigorous justification of the hydrostatic
		approximation}.
	\newblock {\em J. Math. Pures Appl.}, 124:30--58, 2019.

	\bibitem{Li-Titi-Yuan2022}
	Jinkai Li, Edriss S. Titi, and Guozhi Yuan.
	\newblock {The primitive equations approximation of the anisotropic horizontally viscous 3D Navier-Stokes equations}. \newblock {\em J. Differential Equations}, 306:492--524, 2022. 
	
	\bibitem{JLLions1994}
	Jacques-Louis Lions, Roger~Temam, and Shouhong~Wang.
	\newblock {Geostrophic asymptotics of the primitive equations of the
		atmosphere}.
	\newblock {\em Topol. Methods Nonlinear Anal.}, 4(2):253--287, 1994.
	
	\bibitem{Lions1992}
	Jacques-Louis Lions, Roger Temam, and Shouhong Wang.
	\newblock {New formulations of the primitive equations of atmosphere and
		applications}.
	\newblock {\em Nonlinearity}, 5(2):237--288, 1992.
	
	\bibitem{Lions1996}
	Pierre-Louis Lions.
	\newblock {\em {Mathematical Topics in Fluid Mechanics. Volume 1.
			Incompressible Models}}.
	\newblock Oxford Lecture Series in Mathematics and Its Applications, 3. Oxford
	University Press, 1996.
	
	\bibitem{Lions1998}
	Pierre-Louis Lions.
	\newblock {\em {Mathematical Topics in Fluid Mechanics. Volume 2. Compressible
			Models}}.
	\newblock Oxford Lecture Series in Mathematics and its Applications , Vol 2, No
	10. Oxford University Press, 1998.
	
	\bibitem{LT2018a}
	Xin Liu and Edriss~S. Titi.
	\newblock {Local well-posedness of strong solutions to the three-dimensional
		compressible primitive equations}.
	\newblock {\em Arch. Rational Mech. Anal.}, 241:729--764, 2021.
	
	\bibitem{LT2018b}
	Xin Liu and Edriss~S. Titi.
	\newblock {Global existence of weak solutions to the compressible primitive
		equations of atmospheric dynamics with degenerate viscosities}.
	\newblock {\em SIAM J. Math. Anal.}, 51(3):1913--1964, 2019.
	
	\bibitem{LT2018LowMach1}
	Xin Liu and Edriss~S. Titi.
	\newblock {Zero Mach number limit of the compressible primitive equations: Well-prepared initial data}.
	\newblock {\em Arch. Rational Mech. Anal.}, 238:705--747, 2020.
	
	\bibitem{LT2022LowMath2}
	Xin Liu and Edriss~S. Titi.
	\newblock {Zero Mach Number Limit of the Compressible Primitive Equations: Ill-prepared Initial Data}.
	\newblock {\em J. Differential Equations}, 356:1--58, 2023.
	
	\bibitem{Maltese2014}
	David Maltese and Anton{\'{i}}n Novotn{\'{y}}.
	\newblock {Compressible Navier-Stokes equations on thin domains}.
	\newblock {\em J. Math. Fluid Mech.}, 16(3):571--594, 2014.
	
	\bibitem{Matsumura1980}
	Akitaka Matsumura and Takaaki Nishida.
	\newblock {The initial value problem for the equations of motion of viscous and
		heat-conductive gases}.
	\newblock {\em J. Math. Kyoto Univ.}, 20(1):67--104, 1980.
	
	\bibitem{Matsumura1983}
	Akitaka Matsumura and Takaaki Nishida.
	\newblock {Initial boundary value problems for the equations of motion of
		compressible viscous and heat-conductive fluids}.
	\newblock {\em Commun. Math. Phys.}, 89:445--464, 1983.
	
	\bibitem{Petcu2005}
	Madalina~Petcu and Djoko~Wirosoetisno.
	\newblock {Sobolev and Gevrey regularity results for the primitive equations in
		three space dimensions}.
	\newblock {\em Applicable Analysis}, 84(8):769--788, 2005.
	
	\bibitem{Renardy2009}
	Michael Renardy.
	\newblock {Ill-posedness of the hydrostatic Euler and Navier--Stokes equations}.
	\newblock {\em Arch. Rational Mech. Anal.}, 194(3):877--886, 2009.
	
	\bibitem{Richardson1965}
	Lewis~F. Richardson.
	\newblock {\em {Weather Prediction by Numerical Process}}.
	\newblock Dover Publications, Inc.,1965.
	
	\bibitem{Serrin1959}
	James Serrin.
	\newblock {On the uniqueness of compressible fluid motions}.
	\newblock {\em Arch. Rational Mech. Anal.}, 3(1):271--288, 1959.
	
	\bibitem{tang2023derivation}
	Tong Tang and {\v{S}}{\'a}rka Ne{\v{c}}asov{\'a}.
	\newblock {Derivation of the inviscid compressible primitive equations}.
	\newblock {\em Applied Mathematics Letters}, 139:108534, 2023.
	
	\bibitem{Tani1977}
	Atusi Tani.
	\newblock {On the first initial-boundary value problem of compressible viscous
		fluid motion}.
	\newblock {\em Publications of the Research Institute for Mathematical Sciences}, 13(1):193--253, 1977.
	
	\bibitem{Temam1984}
	Roger Temam.
	\newblock {\em {Navier--Stokes Equations: Theory and Numerical Analysis}}.
	\newblock American Mathematical Soc., 2001.
		
	\bibitem{Vasseur2016}
	Alexis~F. Vasseur and Cheng Yu.
	\newblock {Existence of global weak solutions for 3D degenerate compressible
		Navier--Stokes equations}.
	\newblock {\em Invent. Math.}, 206(3):935--974, 2016.
	
	\bibitem{Vaygant1995}
	Vladimir Andreevich Va\v{i}gant and Aleksandr Vasilevic Kazhikhov.
	\newblock {On existence of global solutions to the two-dimensional
		Navier--Stokes equations for a compressible viscous fluid}.
	\newblock {\em Siberian Mathematical Journal}, 36(6):1108--1141, 1995.
	
	\bibitem{Wang2017}
	Fengchao Wang, Changsheng Dou, and Quansen Jiu.
	\newblock {Global existence of weak solutions to 3D compressible primitive equations with
		density-dependent viscosity}.
	\newblock {\em Journal of Mathematical Physics}, 61:021507, 2020.
	
	\bibitem{Washington2005}
	Warren~M. Washington and Claire~L. Parkinson.
	\newblock {\em {An Introduction to Three-Dimensional Climate Modeling}}.
	\newblock University Science Books, 2005.
	
	\bibitem{zpxin1998}
	Zhouping Xin.
	\newblock {Blowup of smooth solutions to the compressible Navier--Stokes
		equation with compact density}.
	\newblock {\em Communications on Pure and Applied Mathematics}, 51(3):229--240, 1998.

	
\end{thebibliography}

\end{document}